\documentclass[11pt]{article}

\usepackage{amsmath}  
\usepackage{amsthm} 
\usepackage{amssymb} 

\usepackage{enumitem}
\usepackage[backref=page]{hyperref}  %

\usepackage{color} 

\newtheorem{thm}{Theorem}
\newtheorem{inspr}[thm]{}

\newenvironment{definitie}{\begin{itemize}\item[ ]\hspace{-26pt}\bf Definition \rm }{\end{itemize}}
\newenvironment{notatie}{\begin{itemize}\item[ ]\hspace{-26pt}\bf Notation \rm }{\end{itemize}}
\newenvironment{voorbeeld}{\begin{itemize}\item[ ]\hspace{-26pt}\bf Example \rm }{\end{itemize}}
\newenvironment{stelling}{\begin{itemize}\item[ ]\hspace{-26pt}\bf Theorem \rm }{\end{itemize}}
\newenvironment{propositie}{\begin{itemize}\item[ ]\hspace{-26pt}\bf Proposition \rm }{\end{itemize}}
\newenvironment{lemma}{\begin{itemize}\item[ ]\hspace{-26pt}\bf Lemma \rm }{\end{itemize}}
\newenvironment{opmerking}{\begin{itemize}\item[ ]\hspace{-26pt}\bf Remark \rm }{\end{itemize}}
\newenvironment{voorwaarde}{\begin{itemize}\item[ ]\hspace{-26pt}\bf Condition \rm }{\end{itemize}}
\newenvironment{probleem}{\begin{itemize}\item[ ]\hspace{-26pt}\bf Problem \rm }{\end{itemize}}
\newenvironment{gevolg}{\begin{itemize}\item[ ]\hspace{-26pt}\bf Corollary \rm }{\end{itemize}}
\newenvironment{niets}{\begin{itemize}\item[ ]\hspace{-26pt}\bf   \rm }{\end{itemize}}

\renewcommand{\Bbb}{\mathbb} 

\newcommand{\defin}{\begin{inspr}\begin{definitie}}  
\newcommand{\edefin}{\end{definitie}\end{inspr}}

\newcommand{\notat}{\begin{inspr}\begin{notatie}}  
\newcommand{\enotat}{\end{notatie}\end{inspr}}

\newcommand{\voorb}{\begin{inspr}\begin{voorbeeld}}  
\newcommand{\evoorb}{\end{voorbeeld}\end{inspr}}

\newcommand{\stel}{\begin{inspr}\begin{stelling}}
\newcommand{\estel}{\end{stelling}\end{inspr}}

\newcommand{\prop}{\begin{inspr}\begin{propositie}}
\newcommand{\eprop}{\end{propositie}\end{inspr}}

\newcommand{\lem}{\begin{inspr}\begin{lemma}}
\newcommand{\elem}{\end{lemma}\end{inspr}}

\newcommand{\opm}{\begin{inspr}\begin{opmerking}}
\newcommand{\eopm}{\end{opmerking}\end{inspr}}

\newcommand{\voorw}{\begin{inspr}\begin{voorwaarde}}
\newcommand{\evoorw}{\end{voorwaarde}\end{inspr}}

\newcommand{\probl}{\begin{inspr}\begin{probleem}}
\newcommand{\eprobl}{\end{probleem}\end{inspr}}

\newcommand{\gev}{\begin{inspr}\begin{gevolg}}
\newcommand{\egev}{\end{gevolg}\end{inspr}}

\newcommand{\nul}{\begin{inspr}\begin{niets}}
\newcommand{\enul}{\end{niets}\end{inspr}}

\newcommand{\bew}{\vspace{-0.3cm}\begin{itemize}\item[ ] \bf Proof\rm: }
\newcommand{\ebew}{\hfill $\qed$ \end{itemize}}

\newcommand{\snl}{\vskip 3pt} 
\newcommand{\ssnl}{\vskip 7pt} 
\newcommand{\nl}{\vskip 12pt} 

\newcommand{\Cal}{\mathcal}
\ssnl

\newcommand{\ot}{\otimes}

\newcommand{\inv}{^{-1}}

\newcommand{\tussenen}{\qquad\quad\text{and}\qquad\quad}

\numberwithin{thm}{section}   
\numberwithin{equation}{section} 

\newcommand{\keepcomment}[1]{}
\newcommand{\oldcomment}[1]{}

\begin{document}

%
%

\centerline{\bf \Large Discrete quantum groups and their duals}
\vspace{13pt}
\centerline{\it A.\ Van Daele \rm ($^{*}$)}
\vspace{20pt}
{\bf Abstract} 
\nl

Discrete quantum groups were introduced as duals of compact quantum groups by Podle\'s and Woronowicz in 1990 \cite{Po-Wo}. Shortly after, they were defined and studied intrinsically in \cite{Ef-Ru} and \cite{VD-discrete}. In 1998, with the introduction of the multiplier Hopf algebras with integrals (also called algebraic quantum groups), the duality between discrete and compact quantum groups became part of the more general duality in the self-dual category of these algebraic quantum groups \cite{VD-alg}. Again a few years later the duality was extended to all locally compact quantum groups (see e.g. \cite{Ku-Va2}).
\ssnl 
In these notes, we give a new and a somewhat updated approach of the theory of discrete quantum groups. In particular, we view them as special cases of algebraic quantum groups. The duality between the compact quantum groups and the discrete quantum groups is seen in this larger context. This has a number of advantages as we will explain. 
\ssnl
On the one hand, we provide quite a bit of information about how all of this fits into the more general theory of algebraic quantum groups and its duality. Occasionally, we even go one step further and look at the most general case of locally compact quantum groups. Also sometimes, we compare with known results in pure Hopf algebra theory. On the other hand however, we have tried to make these notes highly self-contained. 
The aim of these notes in the first place is not to give new results but rather to review known results in a more modern perspective, taking into account recent developments. We believe this may be helpful for people who want to work with compact and discrete quantum groups now. 

\nl
Date: {\it 1 April 2026}
\vskip 7 cm
\hrule

\begin{itemize}[noitemsep]
\item[{($^{*}$)}] Department of Mathematics, KU Leuven, Celestijnenlaan 200B, B-3001 Heverlee, Belgium. E-mail: \texttt{Alfons.VanDaele@kuleuven.be}
\end{itemize}
\newpage

%
%

\setcounter{section}{-1}  

\section{\hspace{-17pt}. Introduction}\label{s:intro} 

The operator algebra approach to quantum groups finds it origin in the attempts to generalize {\it Pontryagin's duality theorem} for abelian locally compact groups (\cite{Po}) to the non-abelian ones. The first results in this direction all describe a dual of a (not necessarily abelian) locally compact group and show how to recover the original group from this dual object. Unfortunately, the nice symmetry, present in the abelian case, had to be given up. With the concept of a {\it Kac algebra}, developed independently by Kac and Vainerman (\cite{Va-Ka}) and Enock and Schwarz (\cite{En-Sw1, En-Sw2}), the duality was again realized within a self-dual category of objects.
\ssnl
However, the work of Drinfel'd (\cite{Dr}) and Jimbo (\cite{Ji}) made people realize that one of the conditions, imposed on the antipode in the theory of Kac algebras, was too restrictive. At the same time, Woronowicz (\cite{Wo1}) constructed the quantum $SU_q(2)$ as a compact quantum group within the operator algebra framework. Also this example did not satisfy the Kac algebra axioms. These two events were the starting point for a new search for a self-dual category of objects, strictly larger than the category of Kac algebras so as to include also the compact quantum groups like Woronowicz' $SU_q(2)$.
\ssnl
The result is now known as the theory of {\it locally compact quantum groups}. The development of the final object is due to Masuda, Nakagami and Woronowicz on the one hand (\cite{Ma-Na, Ma-Na-Wo}) and Kustermans and Vaes on the other hand (\cite{Ku-Va1, Ku-Va2, Ku-Va3}, see also \cite{VD-kerala, VD-sigma}.
\ssnl
An intermediate step towards this result came with the theory of {\it multiplier Hopf $^*$-algebras with positive integrals} (sometimes called $^*$-algebraic quantum groups) (\cite{VD-mha, VD-alg}). It is a category of objects, containing both the discrete and the compact quantum groups and it is also self-dual. It is not general enough to contain all locally compact quantum groups. On the other hand however, it is relatively simple and purely algebraic and the theory is rich enough to incorporate most algebraic features of the general theory, avoiding the more difficult technicalities.
\ssnl
The theory of locally compact quantum groups is now well-developed. It is well-understood and there are interesting and plenty of examples, and there are the common general construction methods. And although there is now also a simplified approach to the theory (\cite{VD-sigma}), it remains technically quite involved and relies heavily on Tomita-Takesaki theory techniques. This is probably one of the reasons why it seems that presently, the theory is not very well absorbed by the mathematical community. Perhaps also the fact that there are so many differences in conventions used by various authors (see the appendix).
\ssnl
Lately there came a renewed interest in various aspects of the theory of discrete and compact quantum groups, but generally speaking, the main references that are used date from the period when the locally compact quantum groups, or even the algebraic quantum groups did not yet existed. On the other hand, there are good reasons to believe that a more modern point of view, where compact and discrete quantum groups are seen as special cases of locally compact quantum groups, in particular as cases of multiplier Hopf $^*$-algebras with positive integrals and where the duality between these two is viewed as a special case of the more general duality, has a lot of advantages. In particular, often the more general approach, certainly within the framework of the algebraic quantum groups, is simpler. But also, one might expect that the recent understanding of locally compact quantum groups can possibly contribute to a better way to treat discrete quantum groups and their duals.
\nl
These considerations take us to the {\it aim} and the {\it content of these notes}.
\ssnl
In the first place, we give a somewhat new and {\it updated approach} of the theory of discrete quantum groups. As mentioned already, discrete quantum groups were first introduced by Podle\'s and Woronowicz  as duals of compact quantum groups (\cite{Po-Wo}). Later, discrete quantum groups were intrinsically defined and studied by Effros and Ruan (\cite{Ef-Ru}) and in \cite{VD-discrete}. Therefore, we now look at the work of Podle\'s and Woronowicz as a way to show that the dual of a compact quantum group (or rather a compact quantum matrix group) is a discrete quantum group. In this sense, their result is a special case of the more general duality between locally compact quantum groups, in particular, between $^*$-algebraic quantum groups. Here, undoubtedly, the general case is simpler than the special case as it is treated in the paper by Podle\'s and Woronowicz. The extra property that the underlying algebra of the dual is a direct sum of matrix algebras follows in the first place from results within the theory of compact quantum groups, but it can also be obtained from the existence of the cointegral in the dual.
\ssnl
It should be mentioned that viewing discrete quantum groups as duals of compact quantum groups is mainly motivated by historical arguments. The notion of a discrete quantum group, obtained as a quantization of a discrete group, following the philosophy of quantization, is a very natural object. And in fact, the development of the theory is easier than that of compact quantum groups. We will comment more on this in the appendix.
\nl
We have made the updating of this approach to discrete quantum groups to a great extend in a self-contained way. However, as we want to emphasize also that the theory is part of a bigger context, we will regularly indicate how the results obtained are special cases of the general results. We will look here in the first place to the theory of algebraic quantum groups. Occasionally, we wil go further and consider the more general case of locally compact quantum groups. Also, when appropriate, we will compare with results in pure Hopf algebra theory.
\ssnl
The reader who has some background in these other and more general theories will benefit from these comments as he will be able to situate the results in a bigger context. The reader with no background of this kind may become anxious to learn more in this direction. He may even profit from these notes as they treat discrete and compact quantum groups more in the spirit of the general theory. However, it is possible to read and profit from this approach if none of this is the case. And we hope, all readers will enjoy the beauty (and simplicity) of this theory. After all, the main intention for distributing these notes is to make this material easily accessible. 

\nl 
\bf Content of the paper \rm 
\nl
In {\it Section} \ref{s:defin}, we start with a more general concept than the one of a discrete quantum group. We take a multiplier Hopf algebra $(A,\Delta)$ and we call it \emph{of discrete type} if there is a left and a right cointegral (see Definitions \ref{defin:1.1a} and \ref{defin:1.2c}). We are particularly interested in the case where a left cointegral $h$ is also a right cointegral. In that case, the element $\Delta(h)$, obtained by applying the coproduct on the left cointegral, plays a very special role. It is a separability idempotent in the multiplier algebra $M(A\ot A)$. It is used to construct a left and a right integral on $A$ and many properties of these integrals are obtained from the properties of the cointegral.
\ssnl
In {\it Section} \ref{s:discr} we apply the results obtained in the first section to the special case of a \emph{discrete quantum group}. This is defined as a multiplier Hopf $^*$-algebra $(A,\Delta)$ of discrete type with the extra constraint that $A$ is a direct sum of full matrix algebras $A_\alpha$ over $\mathbb C$, see Definition \ref{defin:1.2}. In this case, the left cointegral can be constructed from the axioms of a multiplier Hopf $^*$-algebra whose underlying $^*$-algebra is such a direct sum. One of the main results here is that the square $S^2$ of the antipode $S$  leaves the components invariant. It is a homomorphism and has to be implemented by and invertible element on each component $A_\alpha$. As it turns out, it is in fact implemented by a \emph{positive} invertible element $q_\alpha$ in $A_\alpha$. The integrals on $A$ are expressed in terms of the standard traces on the components, using these elements $q_\alpha$. Hence, we see that the integrals and the corresponding data, are all expressible using the data of the separability idempotent $\Delta(h)$. This treatment of discrete quantum groups is somewhat different from the original way as it is done in \cite{VD-discrete}.
\ssnl
In {\it Section} \ref{s:dual}, we consider duality. In the first place, we look at the dual of a discrete quantum group which is a compact quantum group. The duality is seen within the framework of multiplier Hopf algebras with integrals. Then, we show the one-to-one correspondence between representations of the underlying algebra of the discrete quantum group and the corepresentations of the dual compact quantum group. We obtain various properties of compact quantum groups in this way. Observe once more that this approach is different from what is the common practice where first compact quantum groups are studied and then information about the dual is obtained. Our approach in these notes is justified as, after all, we basically study discrete quantum groups. 
\ssnl
In {\it Section} \ref{s:examples}, we simply discuss the classical cases of a discrete group and of the dual of a compact group. This is just to illustrate the preceding treatment of discrete quantum groups. 
In a separate paper \cite{VD-suq2} we give a detailed treatment of the discrete quantum group $su_q(2)$ as another example of the theory of discrete quantum groups as it is studied in this paper. I believe that this is new although the underlying ideas are not. Indeed, in earlier work the discrete quantum group $su_q(2)$ has been studied as the dual of the compact quantum group $SU_q(2)$ as introduced by Woronowicz in \cite{Wo1}, see again \cite{Po-Wo}. In \cite{VD-suq2}  we reverse the order and start with the discrete quantum group. The dual is then the compact quantum group $SU_q(2)$.
\ssnl
In the last section, {\it Section} \ref{s:concl} we  draw some conclusions and discuss possible further research.
\ssnl
Finally, in an {\it Appendix}, we treat some more aspects of the differences in approaches to discrete and compact quantum groups, as well as the different conventions that are used by different authors.
\nl
\bf Notations and conventions \rm
\nl
Except for Section \ref{s:defin}, we only work with $^*$-algebras over the complex numbers.
The identity element in an algebra is denoted by $1$ while the identity in a group will be denoted by $e$. We use $\iota$ to denote the identity map. The linear dual of a vector space $V$ will be denoted by $V'$. 
\ssnl
In most cases, the tensor product is purely algebraic. If this is not the case, we will say so explicitly.
\ssnl
We use the leg-numbering notation for elements in, as well as operators on a tensor product. If e.g.\ $x$ is an element in the tensor product $A\ot A$ of an algebra $A$ (with unit) with itself, then $x_{12}$, $x_{13}$ and $x_{23}$ are the elements in $A\ot A\ot A$ where $x$ is considered as sitting in respectively the first two factors, the first and the third factor and the second and the third factor. At the remaining place, we put the identity $1$ of the algebra. So, we get e.g.\ $x_{12}=x\ot 1$. Similarly for linear operators where at the remaining factor, the identity map $\iota$ is placed.
\ssnl
Because it is a very convenient tool, sometimes in these notes, we will also use the Sweedler notation when $(A,\Delta)$ is a multiplier Hopf algebra. We write
\begin{equation}
\Delta(a)=\sum_{(a)}a_{(1)}\ot a_{(2)}\qquad\text{and}\qquad
(\Delta\ot\iota)\Delta(a)=\sum_{(a)}a_{(1)}\ot a_{(2)} \ot a_{(3)} \label{eqn:1.1}
\end{equation}
when $a$ is an element of $A$. This can be done, provided we multiply in all but one factor with elements from $A$, left or right. The use of the Sweedler notation for multiplier Hopf algebras has been introduced in \cite{Dr-VD}. The technique has been refined in \cite{VD-tools} and in an appendix of \cite{VD-lnalg}. For these discrete quantum groups, it is relatively transparent as one can simply multiply with central projections, in other words, one can restrict to components. 
\nl
\bf Basic references \rm
\nl
The standard references for the theory of Hopf algebras are \cite{Ab} and \cite{Sw}. There is also the more recent work of Radford \cite{Ra1}. 
\ssnl
 The main references for multiplier Hopf algebras and algebraic quantum groups are \cite{VD-mha} and \cite{VD-alg}. Other possible references here are \cite{VD-Zh1} and \cite{VD-part1, VD-part2, VD-part3}, as well as \cite{Ti}. See also the more recent paper \cite{DC-VD} with some new results. For discrete quantum groups we refer to \cite{Po-Wo}], \cite{Ef-Ru} and \cite{VD-discrete} but our treatment here will be closest to \cite{VD-discrete}. For compact quantum groups we have the basic references \cite{Wo2} and \cite{Wo3}, see also \cite{Ma-VD} for an other treatment. 
\ssnl
The main references for the theory of locally compact quantum groups that we use are \cite{Ku-Va2} and \cite{Ku-Va3}. See also the newer and simplified approach in \cite{VD-sigma}.

\nl
\bf Acknowledgments \rm
\nl
I am very grateful to my colleagues and friends, for the warm hospitality while staying at the University of Fukuoka (Japan) where this work was initiated. I am also indebted to the KU Leuven for the opportunity to continue my research after my official retirement.

%

\section{\hspace{-17pt}.  Multiplier Hopf algebras of discrete type}\label{s:defin} 

We start this paper on discrete quantum groups with a more general concept.

\defin\label{defin:1.1a} 
Let $(A,\Delta)$ be a multiplier Hopf algebra. A non-zero element $h$ in $A$ satisfying $ah=\varepsilon(a)h$ for all $a\in A$ is called a \emph{left cointegral}. A non-zero element $k$ in $A$ satisfying $ka=\varepsilon(a)k$ for all $a\in A$ is called a \emph{right cointegral}. 
\edefin

First we have the following easy observation. 

\prop\label{prop:1.2a} 
 If $h$ is a left cointegral satisfying $\varepsilon(h)\neq 0$, then it is also a right cointegral. 
 It can be normalized such that $h^2=h$ and then it is unique.
 \eprop
 
 \bew
 If $h$ satisfies $\varepsilon(h)\neq 0$, we can take a scalar multiple and assume that $\varepsilon(h)=1$. Then $h^2=h$. Moreover, if $k$ is a right cointegral we have $kh=\varepsilon(h)k=\varepsilon(k)h$. So $k=\varepsilon(k)h$. Because $k$ is also non-zero we must have $\varepsilon(k)\neq 0$, we can also scale $k$ so that $\varepsilon(k)=1$ and then $h=k$. 
  \ssnl
 If $h'$ is any other left cointegral, we have $hh'=\varepsilon(h)h'=h'$ and because $h$ is a right cointegral we have $hh'=\varepsilon(h')h$. It follows that $h'=\varepsilon(h')h$ so that $h'$ is a scalar multiple of $h$. Because $h'$ is also assumed to be non-zero, we must have $\varepsilon(h')\neq 0$ and we can scale $h'$ so that $\varepsilon(h')=1$. Then $h'=h$  and we get uniqueness.
 \ebew

It can happen that $\varepsilon(h)=0$ and that $h$ and $k$ are different. Still, also in that case, one can show that cointegrals are unique, see \cite{VD-Zh0}. 
\ssnl
If $(A,\Delta)$ is a \emph{regular} multiplier Hopf algebra and  $h$ is a left cointegral, then $S(h)$ is a right cointegral because for a regular multiplier Hopf algebra, the antipode is a bijective anti-isomorphism of $A$. If $\varepsilon(h)\neq 0$ we will have $S(h)=h$. But also in this case, it can happen that $\varepsilon(h)=0$. In fact that is already possible for finite-dimensional Hopf algebras.
\ssnl
If $(A,\Delta)$ is a multiplier Hopf $^*$-algebra and if $h$ is a left cointegral, we have $h^*h=\varepsilon(h^*)h=\varepsilon(h)h$. If $\varepsilon(h)=0$ we would have $h^*h=0$. In a sense, this means that we have an involution that is not very useful. On the other hand, if $\varepsilon(h)\neq 0$ we see that we can scale $h$ so as to become a self-adjoint idempotent.
Therefore, it makes sense to require that we have $^*$-algebras with the property that $a^*a=0$ implies that $a=0$. This is the case when $A$ is an operator algebra.  We will impose this condition for a discrete quantum group in Section \ref{s:discr}.
\ssnl
Further in this section, we will consider these different cases.

\defin\label{defin:1.2c}
A  multiplier Hopf algebra  is called of \emph{discrete type} if there is a left and a right cointegral.
\edefin

We will now prove several properties of cointegrals. Most of them can be found already in \cite{VD-Zh0}. 

\prop
For a left cointegral $h$ and a right cointegral $k$ we have
\begin{equation*}
\Delta(a)(1\ot h)=a\ot h
\tussenen
(k\ot 1)\Delta(a)=k\ot a
\end{equation*}
for all $a\in A$. 
\eprop

This follows immediately from the definitions.  
\ssnl
Also the following result is still easy to obtain.

\prop\label{prop:2.5a}
Let $(A,\Delta)$ be any multiplier Hopf algebra. Assume that $h$ is a left cointegral in $A$. Then
$(1\ot a)\Delta(h)=(S(a)\ot 1)\Delta(h)$ 
for all $a\in A$. If $k$ is a right cointegral, we have $\Delta(k)(a\ot 1)=\Delta(k)(1\ot S(a))$ for all $a$. 
\eprop

\bew
First remark the following. These equations hold in the multiplier algebra $M(A\ot A)$, except when $(A,\Delta)$ is regular. Then they hold in $A\ot A$.
\ssnl
Using  the Sweedler notation we write
$1\ot a=\sum_{(a)}(S(a_{(1)})\ot 1)\Delta(a_{(2)})$. Then, for all $a\in A$,
\begin{align*}
(1\ot a)\Delta(h)&=\sum_{(a)}(S(a_{(1)})\ot 1) \Delta(a_{(2)}h)\\
&=\sum_{(a)}(S(a_{(1)})\varepsilon(a_{(2)})\ot 1) \Delta(h)\\
&=(S(a)\ot 1)\Delta(h).
\end{align*}
The other equation is proven in a similar way.
\ebew

In the case of a multiplier Hopf $^*$-algebra we can use that $S(a^*)=S^{-1}(a)^*$ to obtain one equation from the other.
\ssnl
Since this is the first time where the Sweedler notation is used, let us explain it a bit more in the following remark. 

\opm
If in all of the above expressions, we multiply in the first factor with any element from the left, we see that we have all these formulas within the algebraic tensor product of $A$ with itself. The various identities we use, like $m(\iota\ot S)\Delta(a))=\varepsilon(a)1$ in  the first equality, multiplication by such an element makes it possible to interpret these formulas in the correct way.
\eopm

\prop\label{prop:2.6a}
Assume that $(A,\Delta)$ is a \emph{regular} multiplier Hopf algebra. If $h$ is a left cointegral then the left and the  right leg of $\Delta(h)$ are  all of $A$. 
\eprop

\bew
Let $h$ be a left cointegral. Consider the subset $A_0$ of $A$ spanned by elements of the form  $(\omega\ot\iota)((c\ot 1)\Delta(h))$ where $c\in A$ and $\omega$ in $A'$. Assume that $\rho$ is a linear functional on $A$ that is $0$ on this subset. Then $(\iota\ot\rho)((c'S(c)\ot 1)\Delta(h))=0$ for all $c,c'$. 
By the previous proposition we get $(\iota\ot\rho)((c'\ot c)\Delta(h))=0$. 
\ssnl
Using the Sweedler notation we see that $\sum_{(h)} c'h_{(1)}\rho(ch_{(2)})=0$. Then we apply the map $ \Delta$ and obtain
\begin{equation*}
\sum_{(h)} \Delta(c')(h_{(1)}\ot h_{(2)})\rho(ch_{(3)})=0.
\end{equation*}
We can multiply with an element $c''$ from the left in the first factor and use that $(A\ot 1)\Delta(A)=A\ot A$. Then we find
\begin{equation*}
\sum_{(h)} ph_{(1)}\ot qh_{(2)}\rho(ch_{(3)})=0
\end{equation*}
for all $p,q,c$. We apply the antipode on the second factor and multiply from the left with an element $r$. Then
\begin{equation*}
\sum_{(h)} ph_{(1)}\ot rS(h_{(2)})S(q)\rho(ch_{(3)})=0.
\end{equation*}
Now $h_{(1)}$ is covered by $p$ and $h_{(2)}$ is covered by $r$. Therefore we can cancel $S(q)$.
We write this as
\begin{equation*}
\sum_{(h)} ph_{(1)}\ot rS(h_{(2)})\ot \rho(\,\cdot\,h_{(3)})=0.
\end{equation*}
Finally we apply the evaluation map on the last two factors and we get for the left hand side
\begin{align*}
\sum_{(h)} ph_{(1)} \rho(rS(h_{(2)})h_{(3)})
=\sum_{(h)} ph_{(1)}\varepsilon(h_{(2)})\rho(r)
=ph\rho(r).
\end{align*}
As this is now equal to $0$ for all $p$ and $r$ we obtain that $\rho=0$. 
Because this is true for all linear functionals that are $0$ on the subspace $A_0$ we must have $A_0=A$. 
\ssnl
ii) Next we claim that also the span of elements $(\omega(\,\cdot\,b)\ot\iota)\Delta(h)$ is all of $A$. Indeed, if $\rho$ is $0$ on all such elements, we must have $(\iota\ot\rho)(\Delta(h)(b\ot 1))=0$ for all $b$. We can multiply with $c$ from the left and cancel $b$ so that again $(\iota\ot\rho)((c\ot 1)\Delta(h))=0$ for all $c$. Then we can proceed as in item i).
\ssnl
iii) 
 In a similar way we can prove that the spaces spanned by elements of the form $(\iota\ot\omega)((1\ot b))\Delta(h)$ and $(\iota\ot\omega)(\Delta(h)(1\ot b))$ are both equal to all of $A$. 
\ebew

Also the two legs of $\Delta(k)$ are all of $A$ for a right cointegral $k$ in a regular multiplier Hopf algebra. 
\ssnl
We can go a little further here, using a standard argument from the general theory of multiplier Hopf algebras.

\prop\label{prop:2.6d}
Let $(A,\Delta)$  be a regular multiplier Hopf algebra with a left cointegral $h$. Assume that $\omega$ is a \emph{faithful} linear functional on $A$. Then any element of $A$ is of the form $(\omega(\,\cdot\,a)\ot \iota)\Delta(h)$ for some $a\in A$.
\eprop

\bew
Let $\rho$ be any linear functional on $A$ satisfying $\rho(c)=0$ for all elements $c$ of the form 
$c=(\omega(\,\cdot\,a)\ot \iota)\Delta(h)$ for some $a\in A$.
Then, for all $a,b$ 
\begin{equation*}
(\omega(\,\cdot\,ab)\ot \rho)(\Delta(h))=0
\end{equation*}
and so $\omega(xb)=0$ where $x=(\iota\ot\rho)(\Delta(h)(a\ot 1))$. Because this holds for all $b$ and as $\omega$ is faithful, it follows that $x=0$. Then $\omega_1(x)=0$ for all $\omega_1$. This means that $\rho$ is actually $0$ on all elements of the form
\begin{equation*}
(\omega_1(\,\cdot\, a)\ot\iota)\Delta(h).
\end{equation*}
By Proposition \ref{prop:2.6a} $A$ is spanned by elements of this form. Hence $\rho$ is $0$ on all of $A$. This proves the result.
\ebew

A similar result will be true with different variations, as well as for right cointegrals.
\ssnl
Here is an important consequence about the structure of the underlying algebra $A$.

\prop
Let $(A,\Delta)$ be a regular multiplier Hopf algebra with a left cointegral $h$. Then $AaA$ is finite-dimensional for all $a$.
\eprop

\bew
i) Consider $c\in A$ and write $\Delta(h)(c\ot 1)=\sum_i p_i\ot q_i$. For all $a$, we have $(1\ot a)\Delta(h)=(S(a)\ot 1)\Delta(h)$ by Proposition \ref{prop:2.5a}. Next let $\omega\in A'$ and put $x=(\omega(\,\cdot\,c)\ot \iota)\Delta(h)$. Then
\begin{align*}
ax
&=(\omega(\,\cdot\,c)\ot \iota)((1\ot a)\Delta(h))\\
&=(\omega(\,\cdot\,c)\ot \iota)((S(a)\ot 1)\Delta(h))\\
&=\omega(S(a)p_i)q_i
\end{align*}
and we see that $ax$ belongs to the finite-dimensional space spanned by the elements $(q_i)$. By Proposition \ref{prop:2.6a} any element of $A$ is spanned by such elements $x$ and therefore we have shown that $aA$ is finite-dimensional for all $a$ in $A$.
\ssnl
ii) We can apply a similar argument with a right cointegral and this gives that $Aa$ is finite-dimensional for all $a$ in $A$. 
\ssnl
iii) Combining the two results we can conclude that $AaA$ is finite-dimensional for all $a\in A$.
\ebew

We can now prove the following characterization property.

\stel\label{stel:2.8a}
Assume that $(A,\Delta)$ is a multiplier Hopf $^*$-algebra of discrete type. If $A$ is an operator algebra then $A$ is a direct sum of finite-dimensional matrix algebras over $\mathbb C$.
\estel

\bew
For all $a$ we have a finite-dimensional two-sided ideal $AaA$. It contains $a$ because $A$ has local units. Because $A$ is an operator algebra, $AaA$ is the direct sum of full matrix algebras. Because any element $a$ belongs to such a sum of matrix algebras, we must have that $A$ is direct sum of matrix algebras. 
\ebew

The next section will be devoted to the study of multiplier Hopf $^*$-algebras of discrete type where the underlying $^*$-algebra is an operator algebra as in this theorem.  They are what we call here  \emph{discrete quantum groups}.
\nl
\bf Invariant integrals on $(A,\Delta)$ \rm
\nl
In what follows, we assume that $(A,\Delta)$ is a  multiplier Hopf algebra of discrete type as in Definition \ref{defin:1.2c}. We also assume now that it is a \emph{regular multiplier Hopf algebra}.

\stel\label{stel:2.4}
Assume that $(A,\Delta)$ is a regular multiplier Hopf algebra with a left cointegral $h$.
Then there exist non-zero linear functionals $\varphi$ and $\psi$ on $A$ satisfying 
\begin{align} (\iota\ot\varphi)\Delta(a)&=\varphi(a)1\\
         (\psi\ot\iota)\Delta(a)&=\psi(a)1
\end{align}
for all $a\in A$. These functionals are faithful.
\estel

\bew
i) We will denote by $A'_r$ the space of linear functionals on $A$ spanned by elements of the form $\omega(c\,\cdot\,)$ where $c\in A$ and $\omega\in A'$. Define $\varphi(x)=\omega(1)$ if $x=(\omega\ot\iota)\Delta(h)$ with $\omega\in A'_r$. As we have seen in the proof of Proposition \ref{prop:2.6a} above all elements in $A$ are of this form. Moreover, if for such an element we have 
$(\omega\ot \iota)\Delta(h)=0$ we must have $\omega=0$ because also the left leg of $\Delta(h)$ is all of $A$.  
\ssnl
Now take any $a\in A$ and write $a=(\omega\ot\iota)\Delta(h)$ with $\omega\in A'_r$. Then for all $\omega'\in A'_r$ we get
\begin{align*}
 \varphi((\omega'\ot\iota)\Delta(a))
           &=\varphi((\omega'\ot\iota)(\omega\ot 1\ot 1)(\Delta \ot\iota)\Delta(h))\\
           &=\varphi((((\omega\ot\omega')\circ\Delta) \ot\iota)\Delta(h))\\
           &=(\omega\ot\omega')\Delta(1)=\omega(1)\omega'(1)\\
           &=\omega'(1)\varphi(a).
\end{align*}
Hence $(\iota\ot\varphi)\Delta(a)=\varphi(a)1$. We also have $\varphi(h)=1$ by definition.
\ssnl
ii) 
The case of the right integral $\psi$ is similar. If we use a right cointegral $k$ to construct it, we get $\psi(k)=1$.
\ssnl
iii) Non-zero invariant linear functionals on a regular multiplier Hopf algebra are automatically faithful.
\ebew

We have used coassociativity of $\Delta$ and the extension of $\Delta$ to the multiplier algebra $M(A)$. All of this is justified because $\omega,\omega'\in A'_r$.

\ssnl
This result, as well as some of the following results, are also found in \cite{VD-Zh0}.

\prop \label{prop:1.12a}
Let $(A,\Delta)$ be a regular multiplier Hopf algebra of discrete type. Let $h$ be a left cointegral and put $k=S\inv(h)$. Let $\varphi$ be the left integral satisfying $\varphi(h)=1$ and $\psi$ the right integral satisfying $\psi(k)=1$. Then $\varphi\circ S=\psi$. If $\varepsilon(h)=1$ we also have $\psi\circ S=\varphi$. 
\eprop
\bew
i) By assumption $S(k)=h$ and $\varphi(S(k))=\varphi(h)=1$. Because $\varphi\circ S$ is a right integral  and takes the same value in $k$ as $\psi$, we will have $\psi=\varphi\circ S$ by the uniqueness of integrals.
\ssnl
ii) By the uniqueness of cointegrals, we have that $S^2(h)$ is a scalar multiple of $h$. Now $\varepsilon(S^2(h))=\varepsilon(h)$. Therefore, if $\varepsilon(h)$ is non-zero, we must have $S^2(h)=h$. Then also S(h)=k and $\varphi=\psi\circ S$ using the argument as in i).
\ebew

We see that the scaling constant $\nu$, defined by $\varphi\circ S^2=\nu \varphi$, has to be equal to $1$ if the cointegral $h$ satisfies $\varepsilon(h)\neq 0$.
\ssnl
If $(A,\Delta)$ is multiplier Hopf $^*$-algebra of discrete type and $A$ is an operator algebra, we will show in the next section that the integrals are positive. 

\nl
\bf The modular automorphisms \rm
\nl
From the general theory of algebraic quantum groups (i.e.\ regular multiplier Hopf algebras with integrals) we know the existence of an automorphism $\sigma$ of $A$ such that $\varphi(ab)=\varphi(b\sigma(a))$ for all $a,b\in A$ when $\varphi$ is a left integral. For the integrals we have obtained in Theorem \ref{stel:2.4} we can obtain the following property of the modular automorphisms.

\prop\label{prop:2.10a}
Let $h$ be a left cointegral and assume that it is also a right cointegral. Then for all $a,b\in A$ we have
$$\varphi(ab)=\varphi(bS^2(a))\qquad \text{and} \qquad \psi(ab)=\psi(bS^{-2}(a)).$$
\eprop

\bew
i) Because $h$ is a left cointegral, we have $(1\ot a)\Delta(h)=(S(a)\ot 1)\Delta(h)$ for all $a$ and because we assume that it is also a right cointegral, we also have $\Delta(h)(a\ot 1)=\Delta(h)(1\ot S(a))$ for all $a$. 
\ssnl
ii) Now let $a\in A$. Then
\begin{align*}
(\iota\ot\varphi)((1\ot a)\Delta(h))
&=(\iota\ot\varphi)((S(a)\ot 1)\Delta(h))\\
&=S(a)=(\iota\ot\varphi)(\Delta(h)(S(a)\ot 1)\\
&=(\iota\ot\varphi)(\Delta(h)(1\ot S^2(a))).
\end{align*}
Because the right leg of $\Delta(h)$ is all of $A$ we get $\varphi(ab)=\varphi(bS^2(a))$ for all $b\in A$.
\ssnl
iii) Similarly we get $\psi(ab)=\psi(bS^{-2}(a))$.
\ebew

\prop
If $h$ is a left cointegral and also a right cointegral we have $\sigma(h)=h$. Then also $S^2(h)=h$.
\eprop
\bew
For all $a\in A$ we have $\varphi(ha)=\varepsilon(a)\varphi(h)=\varphi(ah)$ and therefore $\sigma(h)=h$. Because $S^2$ and $\sigma$ coincide we also get $S^2(h)=h$.
\ebew

The property $S^2(h)=h$ is true for any left cointegral if $\varepsilon(h)\neq 0$. In that case, $h$  is also a right cointegral. What we obtain here is that $S^2(h)=h$ holds if $h$ is also a right cointegral. It is not clear if this would imply already that $\varepsilon(h)\neq 0$. 
\ssnl
We also have a converse result.

\prop\label{prop:2.9a}
Assume that $(A,\Delta)$ is a regular multiplier Hopf algebra with a left cointegral $h$. Let $\sigma$ be the modular automorphism of the left integral $\varphi$. If $\sigma(h)=h$ we must have that $h$ is also a right cointegral.
\eprop

\bew
If $\sigma(h)=h$ we have $\varphi(ha)=\varphi(ah)=\varepsilon(a)\varphi(h)=\varepsilon(a)$. This means that 
$$\varphi(hab)=\varepsilon(ab)=\varepsilon(a)\varepsilon(b)=\varepsilon(a)\varphi(hb).$$ From the faithfulness of $\varphi$ we find that $ha=\varepsilon(a)h$.
\ebew

\bf The modular element \rm
\nl
Again from the general theory, we know the existence of a multiplier $\delta\in M(A)$ given by $\psi(a)=\varphi(a\delta)$ for all $a$ (assuming that the $\psi$ and $\varphi$ are normalized in such a way that $\psi=\varphi\circ S$). We also have $(\varphi\ot\iota)\Delta(a)=\varphi(a)\delta$ for all $a$. In particular, because we use the normalization $\varphi(h)=1$, we get here $\delta=(\varphi\ot\iota)\Delta(h)$. Again we can not say much more except when, as in Proposition \ref{prop:2.9a} , a left cointegral is also a right cointegral.

\prop\label{prop:1.14c}
If $h$ is also a right cointegral, then for all $a$ we have $S^4(a)=\delta\inv a\delta$.
\eprop

\bew
For $a\in A$ we have
\begin{align*}
\delta a
&=(\varphi\ot\iota)((1\ot a)\Delta(h))\\
&=(\varphi\ot\iota)((S(a)\ot 1)\Delta(h))\\
&=(\varphi\ot\iota)(\Delta(h)(S^3(a)\ot 1)\\
&=(\varphi\ot\iota)(\Delta(h)(1\ot S^4(a))\\
&=\delta S^4(a).
\end{align*}
\ebew

This is Radford's formula proven in \cite{Ra0} for the fourth power of the antipode because $\widehat \delta=1$. That $\widehat \delta=1$ follows from the general theory.  Indeed, we know that $\langle a,\widehat\delta\rangle=\varepsilon(\sigma\inv(a))$, see e.g.\ Proposition 2.3 in \cite{VD-part1} or Proposition 1.2.3 in \cite{VD-lnalg}. In this case we have $S^2=\sigma$ and because $\varepsilon(S^{-2}(a))=\varepsilon(a)$ we get $\langle a,\widehat\delta\rangle=\varepsilon(a)$ and therefore $\widehat\delta=1$.
\ssnl

\nl
\bf The element $\Delta(h)$ as a separability idempotent \rm
\nl 
In this item, we again assume that $h$ is a left cointegral and \emph{that it satisfies $\varepsilon(h)=1$}. Then it also a right cointegral and we can assume that it is an idempotent as we saw in Proposition \ref{prop:1.2a}. It follows that $\Delta(h)$ is an idempotent in $M(A\ot A)$. It plays and important role in the next section on discrete quantum groups. 

\prop
The element $\Delta(h)$ is a (regular) separability idempotent in $M(A\ot A)$ in the sense of \cite{VD-si1}.
\eprop

\bew
i) First we look at the setting in \cite{VD-si1}. For the algebras $B$ and $C$ we take $A$ and we let $E=\Delta(h)$. The algebra $A$ is non-degenerate and $\Delta(h)$ is an idempotent in $M(A\ot A)$. Because $A$ is a regular multiplier Hopf algebra we know that elements of the form $\Delta(h)(1\ot a)$ and $\Delta(h)(a\ot 1)$ belong to $A\ot A$  for all $a\in A$ and similarly with $a$ on the other side. This shows that Assumption 1.1 of \cite{VD-si1} is fulfilled and that $\Delta(h)$ is regular in the sense of that assumption.
\ssnl
ii) We have shown that the legs of $\Delta(h)$ are all of $A$ and hence it is full. Also Proposition 1.4 in \cite{VD-si1} here says that $\Delta(h)(1\ot a)=0$ only if $a=0$, and similarly $(a\ot 1)\Delta(h)=0$ only if $a=0$. This we knew already because we have a multiplier Hopf algebra and $h\neq 0$.
\ssnl
iii) Finally we come to Definition 1.5 of \cite{VD-si1}. Since we know here already that 
\begin{equation*}
\Delta(h)(a\ot 1)=\Delta(h(1\ot S(a)) \tussenen (1\ot a)\Delta(h)=(S(a)\ot 1)\Delta(h)
\end{equation*}
for all $a$, and because $S$ is bijective, we finally get that $\Delta(h)$ satisfies all the requirements and so it is a regular separability idempotent.
\ebew

The maps $S:B\to C$ and $S':C\to B$ we have in Proposition 1.7 of \cite{VD-si1} are precisely the antipode $S$ of $A$. If we look at the Proposition 1.8 in \cite{VD-si1}, we get here
\begin{equation*}
m(S\ot \iota)(\Delta(h)(1\ot a))
=\sum_{(h)} S(h_{(1)})h_{(2)}a
=\varepsilon(h)a
\end{equation*}
and consequently $m(S\ot \iota)(\Delta(h)(1\ot a))=a$ for all $a$.
\ssnl
Similarly
\begin{equation*}
m(\iota\ot S)((a\ot 1)\Delta(h))
=\sum_{(h)} ah_{(1)})S(h_{(2)})
=\varepsilon(h)a
\end{equation*}
and consequently $m(\iota\ot S)((\ot\iota)\Delta(h))=a$ for all $a$.

\opm
In the case where the left cointegral $h$ satisfies $\varepsilon(h)=1$, so that $\Delta(h)$ is a separability idempotent, we can apply general results from \cite{VD-si1} to obtain the existence of $\varphi$ and $\psi$, as in Theorem \ref{stel:2.4} above. However, we only would get $(\iota\ot\varphi)\Delta(h)=1$ and $(\psi\ot\iota)\Delta(h)=1$ and we would still need the argument given in the proof of Theorem \ref{stel:2.4} to obtain left and right invariance on all of $A$.
\ssnl
Also the results in Proposition \ref{prop:2.10a} and Proposition \ref{prop:1.14c} already follow from the general theory in \cite{VD-si1}.
\eopm

In the next section, we will give an explicit formula for $\Delta(h)$ and refer to this item on separability idempotents.

\section{\hspace{-17pt}. Discrete quantum groups}\label{s:discr}

In this paper, a discrete quantum group is viewed as a \emph{discrete quantum space with a comultiplication}. Therefore we first give the definition of a discrete quantum space as we view it here.

\nl
\bf Discrete quantum spaces \rm
\nl
Let $I$ be an index set. Associate to each member $\alpha$ of $I$ a natural number $n(\alpha)$ and the $^*$-algebra $A_\alpha$ of $n(\alpha)\times n(\alpha)$ complex matrices (with the usual $^*$-algebra structure). Let $A$ be the (algebraic) direct sum of these algebras. Elements of $A$ are written as $(a_\alpha)_{\alpha\in I}$ where $a_\alpha\in A_\alpha$ for each $\alpha\in I$. By assumption only finitely many components $a_\alpha$ are non-zero. This direct sum $A$ is again a $^*$-algebra in the obvious sense. Observe that the algebra has no unit, except when $I$ is a finite set, but that the product is non-degenerate (as a bilinear form). We say that the algebra $A$ is {\it a direct sum of matrix algebras}. The algebras $A_\alpha$ are considered as sitting inside $A$. The common notation is $\sum_{\alpha\in I} \oplus A_\alpha$.
\ssnl
Elements in the multiplier algebra $M(A)$ of $A$ can be written as $(a_\alpha)_{\alpha\in I}$ where $a_\alpha\in A_\alpha$ for each $\alpha\in I$, but without further restrictions. This algebra is called \emph{the direct product of the algebras $A_\alpha$}. Here the usual notation is $\prod_{\alpha\in I} A_\alpha$.
\ssnl
We also consider the tensor product $A\ot A$ of $A$ with itself and the multiplier algebra $M(A\ot A)$ of $A\ot A$. The algebra $A\ot A$ is the direct sum of the tensor product algebras $A_\alpha\ot A_\beta$ where $\alpha,\beta\in I$ and of course the multiplier algebra $M(A\ot A)$ is the direct product of all these algebras. So elements in $M(A\ot A)$ are written as $(a_{\alpha\beta})_{\alpha,\beta\in I}$ with $a_{\alpha\beta}\in A_\alpha\ot A_\beta$. When only finitely many components are non-zero, we get the elements of $A\ot A$.
\ssnl
The identity in $A_\alpha$ will be denoted by $1_\alpha$ and it will also be considered as a minimal central projection in $A$.
\nl
Now we come to quantizing the 'product structure'. So we introduce the notion of a coproduct on such an algebra. 
\nl 
\bf Discrete quantum groups \rm
\nl

\defin\label{defin:1.1}
Let $A$ be a direct sum of matrix algebras as above. A {\it coproduct} (or comultiplication) on $A$ is a $^*$-homomorphism $\Delta: A \to M(A\ot A)$ such that
\ssnl
i)\ $\Delta(A)(1\ot A)\subseteq A\ot A$ and $(A\ot 1)\Delta(A)\subseteq A\ot A$,\newline
ii)\ $(\iota\ot\Delta)\Delta(a)=(\Delta\ot\iota)\Delta(a)$ for all $a\in A$.
\edefin
 Recall that we use $\iota$ to denote the identity map.
\ssnl
The condition i) is used to give a meaning to the condition in ii). Indeed, if we multiply the left hand side of the equation on the left with an element of the form $c\ot 1\ot 1$ where $c\in A$, the expression is $(\iota\ot\Delta)((c\ot 1)\Delta(a))$ and this is well-defined in  $M(A\ot A\ot A)$. Similarly, we multiply the right hand side of the equation from the right with an element of the form $1\ot 1\ot b$ with $b\in A$ and then the expression is $(\Delta\ot\iota)(\Delta(a)(1\ot b))$. This is how the assumptions in i) are used to interpret the equation in ii) in the multiplier algebra $M(A\ot A\ot A)$. 
\ssnl
Condition ii) is called {\it coassociativity}.
\ssnl
The above definition of a coproduct is also possible when $A$ is any $^*$-algebra with a non-degenerate product (see \cite{VD-mha}). If we do not assume the conditions in i), but only assume that $\Delta$ is non-degenerate (i.e.\ that $\Delta(A)(A\ot A)=A\ot A$), it is still possible to give a meaning to coassociativity, see e.g.\ the note on coassociativity \cite{VD-coass, VD-coass1}. This is not so hard to see in the case of a direct sum. However, as we will see, the  conditions in i) are quite natural and similar conditions are often taken as part of the notion of a coproduct itself.
\nl
Then, the following is the main definition for our approach here.

\defin\label{defin:1.2}
A {\it discrete quantum group} is a pair $(A,\Delta)$ of a $^*$-algebra $A$ that is a direct sum of matrix algebras and a coproduct $\Delta$ on A (as in Definition \ref{defin:1.1}) such that $(A,\Delta)$ is a multiplier Hopf $^*$-algebra.
\edefin

We know that there is a unique {\it counit} and a unique {\it antipode} $S$. The counit is defined as a linear map $\varepsilon:A\to \Bbb C$ satisfying 
\begin{equation*}
(\varepsilon\ot\iota)\Delta(a)=a \qquad \text{and}\qquad 
         (\iota\ot \varepsilon)\Delta(a)=a
\end{equation*}
for all $a$ in $A$. The {\it antipode} is defined as a linear map $S: A \to A$ such that 
$$m(S\ot\iota)\Delta(a)=\varepsilon(a)1 \qquad \text{and}\qquad  m(\iota\ot S)\Delta(a)=\varepsilon(a)1$$
for all $a$ in $A$ (where $m$ is the multiplication map, seen as a linear map from $A\ot A$ to $A$). 
\ssnl
Also these formulas make sense in the multiplier algebra $M(A)$. One just has to multiply, left or right (depending on the case) with an element in $A$.
\ssnl
From the fact that the comultiplication is a homomorphism, it follows that $\varepsilon$ is a homomorphism and that $S$ is a anti-homomorphism. Because $\Delta$ is a $^*$-homomorphism, it follows that $\varepsilon$ is also a $^*$-homomorphism and that $S$ satisfies the equation $S(S(a)^*)^*=a$ for all $a$ in $A$. In particular, $S$ is bijective and $S(a^*)=S^{-1}(a)^*$ for all $a$. We see that $S$ is a $^*$-map if and only if $S^2=\iota$.
\ssnl

Before we continue, we give the \emph{basic example}, motivating the main definition. It shows that the counit is indeed the dual of the unit element in a group and that the antipode is what corresponds to the inverse in a group. The example justifies Definition \ref{defin:1.2} as a natural definition of a discrete quantum group. 

\voorb\label{voorb:2.3a}
Let $G$ be a (discrete) group. Consider the algebra of complex functions on $G$ with finite support (with pointwise operations). Define a coproduct $\Delta$ on $A$ by $\Delta(f)(p,q)=f(pq)$ where $p,q\in G$. Then $(A,\Delta)$ is a discrete quantum group.
\evoorb

In this case, the index set is equal to $G$ and $n(p)=1$ so that $A_p=\Bbb C 1$ for all $p\in G$.
\ssnl
In Section \ref{s:examples} we also consider some other important special cases and we refer to \cite{VD-suq2} for a more specific example.
\nl
There is a lot more to say about this definition. This is done in the appendix. But we give here some comments already. 

\opm
i) If $(A,\Delta)$ is a pair of a $^*$-algebra with a coproduct as in Definition \ref{defin:1.2}, we can look at the linear maps $T_1$ and $T_2$ defined on $A\ot A$ by 
$$ T_1(a\ot b)=\Delta(a)(1\ot b) \qquad \text{and} \qquad T_2(a\ot b)=(a\ot 1)\Delta(b).$$ 
By the assumption on the coproduct, they map into $A\ot A$. It can be shown however that they are bijective when a counit and an antipode exists. In that case, the inverses are given by
\begin{align} T_1^{-1}(a\ot b)&=((\iota\ot S)\Delta(a))(1\ot b) \\       
         T_2^{-1}(a\ot b)&=(a\ot 1)((S\ot\iota) \Delta(b)).
\end{align} 
So, the pair $(A,\Delta)$ is a multiplier Hopf $^*$-algebra as defined in \cite{VD-mha}.
\ssnl
ii) Conversely, if $(A,\Delta)$ is a multiplier Hopf $^*$-algebra and if $A$ is a direct sum of matrix algebras, then it is a discrete quantum group in the sense of Definition \ref{defin:1.2}. The reason is that in any multiplier Hopf $^*$-algebra, there exists a counit and an antipode (see \cite{VD-mha}.
\eopm

So, a discrete quantum group is nothing else but a multiplier Hopf $^*$-algebra where the underlying algebra is a direct sum of matrix algebras (cf.\ the definition given of a discrete quantum group in \cite{VD-discrete}).
\ssnl
In the previous section, we developed the theory of multiplier Hopf algebras of discrete type. In Theorem \ref{stel:2.8a}, we have conditions to assure that the underlying algebra is a direct sum of matrix algebras. In other words, we get a discrete quantum group (as in definition \ref{defin:1.1}). Next we will see that also conversely, if we have a discrete quantum group, it has cointegrals and so it is a multiplier Hopf algebra of discrete type (as in Definition \ref{defin:1.2c}.
\nl
\bf The cointegral in a discrete quantum group\rm
\nl
The following result is fairly obvious and we only include a proof for completeness.

\prop\label{prop:2.5b}
Assume that $A$ is a direct sum of matrix algebras as considered before and let $\varepsilon$ be a $^*$-homomorphism from $A$ to $\mathbb C$. Then there exists a unique non-zero self-adjoint idempotent $h$ in $A$ with the property that $ah=\varepsilon(a)h$ for all $a\in A$. We also have $ha=\varepsilon(a)h$ for all $a\in A$.
\eprop

\bew
The kernel $A_0$ of $\varepsilon$ is a two-sided, $^*$-invariant ideal of $A$ with codimension $1$. As $A$ is a direct sum of matrix algebras, we must have a central projection $p\in M(A)$ such that $A_0=Ap$. When $h=1-p$, then $h\in A$ and $\varepsilon(h)=1$ and so $\varepsilon(ah-\varepsilon(a)h)=0$. This means that $ah-\varepsilon(a)h$ belongs to the kernel $Ap$ for all $a$. However, because $hp=0$, also $ah-\varepsilon(a)h=0$ and so $ah=\varepsilon(a)h$.  By taking adjoints, we get also $ha=\varepsilon(a)h$ for all $a$. If now $k$ is another such element, then $hk=\varepsilon(h)k=k$ and also $hk=\varepsilon(k)h=h$ and therefore $h=k$. 
\ebew

Remark that in the proof above, we have applied $\varepsilon$ on $M(A)$. This is no problem as we can easily extend $\varepsilon$ to a $^*$-homomorphism on $M(A)$ by $\varepsilon(m)=\varepsilon(am)$ where $a$ is in $A$ and satisfies $\varepsilon(a)=1$.
\ssnl
If $(A,\Delta)$ is a discrete quantum group as in Definition \ref{defin:1.2}, we not only have that $h^*=h$ and $h^2=h$, but also $S(h)=h$ because $\varepsilon(S(a))=\varepsilon(a)$ for all $a$.
\ssnl
A trivial consequence is the following.

\prop
When $(A,\Delta)$ is a discrete quantum group, then it is a multiplier Hopf algebra of discrete type. The left cointegral $h$ satisfies $\varepsilon(h)=1$, it is an idempotent and $S(h)=h$. It is also a right cointegral.
\eprop

It turns out to give the integral on the dual - see Proposition \ref{prop:3.3} in the next section.
\ssnl

We can apply all the results of the previous section because $\varepsilon(h)\neq 0$.

\notat
i) The support of the counit $\varepsilon$ is a distinct one-dimensional component. We will denote the corresponding index with $e$, reminding us of the identity element in a group. 
\ssnl
ii) Let $\alpha$ be any index and consider the corresponding component $A_\alpha$. Because $S$ is an anti-isomorphism we have an index $\overline\alpha$ such that $S(A_\alpha)=A_{\overline\alpha}$. 
\enotat

For $a\in A_\alpha$ we have $S(S(a_\alpha)^*)=a^*$.  Because all the components are $^*$-invariant, this implies  $S(A_{\overline\alpha})=A_\alpha$. In particular, $S^2(A_\alpha)=A_\alpha$ for all $\alpha$. We will study these maps later in greater detail. 
\ssnl
First we have the following result about the element $\Delta(h)$ where $h$ is the cointegral given in Proposition \ref{prop:2.5b}.
Recall that we use $1_\alpha$ for the identity of $A_\alpha$ and consider it as a central element in $A$.
\prop\label{prop:2.8a}
We have 
\begin{equation*}
\Delta(h)(1\ot 1_\alpha)=\Delta(h)(1_{\overline\alpha}\ot 1)=\Delta(h)(1_{\overline\alpha}\ot 1_\alpha).
\end{equation*}
\eprop

This follows from the property that 
$S(1_\alpha)=1_{\overline\alpha}$ 
for all $\alpha$.

\nl
\bf The map  $a\mapsto S^2(a)$ and the scaling group \rm
\nl
We  consider the implementation of the square $S^2$ of the antipode $S$ in the following proposition.

\prop\label{prop:2.11a}
Let $\omega$ be any trace on $A$. Define $q=(\omega\ot\iota)\Delta(h)$.  The element $q$ belongs to $M(A)$ and it satisfies  $aq=qS^2(a)$ for all $a\in A$. If $\omega$ is positive, then $q$ is positive and if $\omega$ is faithful, then $q$ is invertible in $M(A)$. 
\eprop
\bew
i) Using the formulas in Proposition \ref{prop:2.5a} and the fact that $\omega$ is a trace, we find
\begin{align*} aq&=(\omega\ot\iota)((1\ot a)\Delta(h))\\
                &=(\omega\ot\iota)((S(a)\ot 1)\Delta(h))\\
                &=(\omega\ot\iota)(\Delta(h)(S(a)\ot 1))\\
                &=(\omega\ot\iota)(\Delta(h)(1\ot S^2(a)))\\
                &=qS^2(a)
\end{align*}
for all $a$. 
\ssnl
ii) As before, we denote by $q_\alpha$ the $\alpha$-component of $q$. When $\omega$ is positive, this will be a positive element in $A_\alpha$. When $\omega$ is faithful, all these components will be non-zero. As they implement $S^2$, they are invertible. Then also $q$ has an inverse in $M(A)$. 
\ebew

\opm
i) There are plenty of traces on $A$, but on each component, the trace is unique up to a scalar. This implies that the components $q_\alpha$ of the element $q$ above are all unique up to a scalar. That is also a consequence of the fact that these components implement the isomorphisms $S^2$, restricted to the components.
\ssnl
ii) Any isomorphism of the algebra of $n\times n$-matrix algebra is implemented by an invertible matrix. Because here $S^2(a)^*=S^{-2}(a^*)$ it is easy to see that it can be implemented by a self-adjoint matrix. That actually we get an implementation by a positive matrix is not automatic but appears to be true in this situation. The reason is that already $S$ satisfies $S(a)^*=S\inv(a^*)$.
\eopm

In a similar way we can consider $p=(\iota\ot\omega)\Delta(h)$ for any trace. Now we have $ap=pS^{-2}(a)$ for all $a\in A$. Again it will be a positive invertible element $p\in M(A)$ if the trace is positive and faithful. Because also $aq\inv=q\inv S^{-2}(a)$ (where $q$ is the element obtained in the previous proposition), we must have that each component $p_\alpha$ of $p$ is a scalar multiple of the corresponding component $q_\alpha$ of $q$.
\nl
We will now look for the scaling group. Recall that algebraic quantum groups admit an \emph{analytic structure}, see \cite{Ku}.
In particular, there is a one-parameter group $(\tau_t)$ of automorphisms of $A$ with the property that all elements in $A$ are analytic. It is characterized by the fact that $S^2=\tau_{-i}$. This is called the \emph{scaling group}. 
\ssnl
In the case of a discrete quantum group, we get the following result.

\prop \label{prop:2.12a}
Let $q$ be a positive invertible element in $M(A)$ satisfying $aq=qS^2(a)$ for all $a\in A$.
Then  $\tau_t(a)=q^{-it} a q^{it}$ for all $t$ and $a\in A$.
\eprop
\bew
Because $S^2(a)=q\inv a q$ and $\tau_{-i}=S^2$ we must have $\tau_t(a)=q^{-it} a q^{it}$.
\ebew

The scaling group leaves all components of $A$ globally invariant and $\tau_t(a)=q_\alpha^{-it} a q_\alpha^{it}$ for all $a\in A_\alpha$. Because $q_\alpha$ is unique up to a scalar, the expression on the right does not depend on the choice of $q$ (as expected).

\nl
\bf The unitary antipode $R$ \rm
\nl
The \emph{unitary antipode} $R$ is characterized so that  $S=R\tau_{-\frac{i}2}$. This formula is thought of as the \emph{polar decomposition of the antipode}. These results are also found in \cite{Ku-VD} and in the more recent survey paper \cite{VD-part3}. The polar decomposition of the antipode is also a feature of the more general locally compact groups, see the references given in the corresponding item in the introduction.
\ssnl
For the unitary antipode of a discrete quantum group,  we get the following results. Recall that for each index $\alpha$ there is an index $\overline\alpha$ such that $S:A_\alpha\to A_{\overline\alpha}$.

\lem 
Take any vector $\xi_0$ in $\Cal H_\alpha \ot \Cal H_{\overline\alpha}$. Define a conjugate linear map $x:\Cal H_\alpha \to \Cal H_{\overline\alpha}$ by
\begin{equation}
\langle x\xi,\eta\rangle = \langle\Delta(h)\xi_0,\xi\ot\eta \rangle.\label{eqn:2.3}
\end{equation}
Then $xa=S(a)^*x$  and $x^*xa=S^2(a)x^*x$ for all $a\in A_\alpha$.
\elem

\bew
i) Recall that $A_\alpha$ acts on the finite-dimensional space $\mathcal H_\alpha$ while $A_{\overline\alpha}$ acts on $\mathcal H_{\overline\alpha}$. Then $\Delta(h)\xi_0$ again belongs to $\Cal H_\alpha \ot \Cal H_{\overline\alpha}$. And because these spaces are finite-dimensional, the formula (\ref{eqn:2.3}) defines a conjugate linear map from $\mathcal H_\alpha$ to $\mathcal H_{\overline\alpha}.$
\ssnl
ii) Now for $a\in A_\alpha$ we get
\begin{align*}
\langle xa\xi,\eta\rangle
&=\langle \Delta(h)\xi_0,a\xi\ot\eta\rangle\\
&=\langle (a^*\ot 1)\Delta(h)\xi_0,\xi\ot\eta\rangle\\
&=\langle (1\ot S\inv(a^*))\Delta(h)\xi_0,\xi\ot\eta\rangle\\
&=\langle (1\ot S(a)^*)\Delta(h)\xi_0,\xi\ot\eta\rangle\\
&=\langle \Delta(h)\xi_0,\xi\ot S(a)\eta\rangle\\
&=\langle x\xi,S(a)\eta\rangle=\langle S(a)^*x\xi,\eta\rangle
\end{align*}
and we see that $xa=S(a)^*x$.
\ssnl
iii) If we take adjoints we get 
$x^*S(a)=a^*x^*$. With $b=S(a)$ we get 
$$x^*b=S\inv(b)^*=S(b^*).$$ 
Now with $b=S(c)^*$ we find $x^*S(c)^*=S(S(c))x^*$. Combining this with the original formula we get 
\begin{equation*}
x^*xa=x^*S(a)^*x=S^2(a)x^*x.
\end{equation*}
\vspace*{-20pt}
\ebew

If $q$ is a positive invertible element satisfying $aq=qS^2(a)$ we will have that $x^*x$ is a scalar multiple of $q_\alpha\inv$. In particular, given such an element $q$ we can find a vector $\xi$ such that  the corresponding $x$  satisfies $x^*x=q_\alpha\inv$.
\ssnl
Then consider the polar decomposition of $x$. It is of the form $x=J_\alpha q_\alpha^{-\frac12}$ where $J_\alpha$ is a conjugate linear isomorphism of $\mathcal H_\alpha$ to $\mathcal H_{\overline\alpha}$.

\prop
The unitary anitpode $R$ maps $A_\alpha$ to $A_{\overline\alpha}$ and for $a\in A_\alpha$ we have $R(a)=J_\alpha(a^*)J_\alpha$.
\eprop
\bew
For $a\in A_\alpha$ we have
\begin{equation*}
S(a)^*=xax\inv=J_\alpha q_\alpha^{-\frac12} a q_\alpha^\frac12 J_\alpha.
\end{equation*}
so that $S(a)=J_\alpha (q_\alpha^{-\frac12} a q_\alpha^\frac12)^* J_\alpha$.
On the other hand, by definition 
\begin{equation*}
S(a)=R(\tau_{-\frac{i}2}(a))=R( q_\alpha^{-\frac12} a q_\alpha^\frac12).
\end{equation*}
It follows that $R(c)=J_\alpha c^* J_\alpha$ for all $c\in A_\alpha$.
\ebew
\nl
\bf A formula for the components of $\Delta(h)$\rm
\nl
We consider an element $q$, positive and invertible in $M(A)$ satisfying $aq=qS^2(a)$ for all $a\in A$ as before. In this item we fix an index $\alpha$ in $I$ and the component $q_\alpha$. It is a positive invertible matrix in $A_\alpha$ satisfying $aq_\alpha=q_\alpha S^2(a)$ for all $a\in A_\alpha$. 
\ssnl
Now we choose matrix elements $(e_{ij})$ in $A_\alpha$ and we assume that $q_\alpha=\sum_{i} \lambda_i e_{ii}$. This is possible because $q_\alpha$ is a positive element. All the eigenvalues $\lambda_i$ are strictly positive as $q_\alpha$ is invertible.

\prop\label{prop:2.16a}
We have 
\begin{equation*}
\Delta(h)(1_{\overline\alpha}\ot 1_\alpha)=\frac1{c}\sum_{ij} \lambda_i\,S(e_{ij})\ot e_{ji}	
\end{equation*}	
where $c=\sum_i \lambda_i$.
\eprop

\bew
Define $H$ in $A_{\overline\alpha}\ot A_\alpha$ by
\begin{equation*}
H=\sum_{i,j} \lambda_i\,S(e_{ij})\ot e_{ji}.
\end{equation*}
\ssnl
i) First we show that $H(a\ot 1)=H(1\ot S(a))$ for all $a$. For any matrix element $e_{k\ell}$ we find
\begin{align*}
H(1\ot e_{k\ell})
&=\sum_{i,j} \lambda_i\,S(e_{ij})\ot e_{ji}e_{k\ell}\\
&=\sum_{j} \lambda_k S(e_{kj})\ot e_{j\ell}\\
H(S\inv(e_{k\ell})\ot 1)
&=\sum_{i,j} \lambda_i\,S(e_{ij})S\inv(e_{k\ell})\ot e_{ji}\\
&=\sum_{i,j} \lambda_i \lambda_\ell\inv\lambda_k S(e_{ij})S(e_{k\ell})\ot e_{ji}\\
&=\sum_{i,j} \lambda_i\lambda_\ell\inv\lambda_k S(e_{k\ell}e_{ij})\ot e_{ji}\\
&=\sum_{j} \lambda_\ell\lambda_\ell\inv\lambda_k S(e_{kj})\ot e_{j\ell}\\
&=\sum_{j} \lambda_k S(e_{kj})\ot e_{j\ell}.
\end{align*}
We have used that $aq_\alpha=q_\alpha S^2(a)$ and that $q_\alpha=\sum_i \lambda_i e_{ii}$. We see that indeed 
$$H(1\ot e_{k\ell})=H(S\inv(e_{k\ell})\ot 1).$$
Then $H(S\inv(a)\ot 1)=H(1\ot a)$ for all $a\in A_\alpha$ and so also $H(a\ot 1)=H(1\ot S(a))$ for all $a\in A_{\overline\alpha}$.
\ssnl
ii) Then we get
\begin{align*}
H\Delta(h)
&=\sum_{(h)} H(h_{(1)}\ot h_{(2)})\\
&=\sum_{(h)}H(1\ot  S(h_{(1)})h_{(2)})\\
&=\varepsilon(h)H=H.
\end{align*}
On the other hand
\begin{align*}
H\Delta(h)
&=\sum_{i,j} \lambda_i (S(e_{ij})\ot e_{ji})\Delta(h)\\ 
&=\sum_{i,j} \lambda_i(1\ot e_{ji}e_{ij})\Delta(h)\\
&=\sum_{i,j} \lambda_i (1\ot e_{jj})\Delta(h)\\
&=c\Delta(h)(1\ot 1_\alpha).
\end{align*}
where $c=\sum_i \lambda_i$.
\ssnl
iii) It follows that $H=c\Delta(h)(1_{\overline\alpha}\ot 1_\alpha)$. 
\ebew

It is natural that this constant is present because we have only used that $q_\alpha$ implements $S^2$ and that is then true for any scalar multiple of $q_\alpha$. Further observe that indeed, the right hand side does not depend on the choice of $q_\alpha$ as it should.

\notat\label{notat:2.15a}
In what follows we will use $\omega_\alpha$ for the standard trace on $A_\alpha$, defined by $\omega_\alpha(1_\alpha)=n_\alpha$ when $n_\alpha$ is the dimension of the matrix algebra $A_\alpha$. 
\enotat
Because $S$ is an anti-isomorphism from $A_\alpha$ to $A_(\overline\alpha)$ we will have that $\omega_{\overline\alpha}\circ S =\omega_\alpha$.
\ssnl
We get the following consequence of Proposition \ref{prop:2.16a}.

\prop
Each component of $\Delta(h)$ is one-dimensional.
\eprop

\bew
Then
\begin{align*}
(\omega_{\overline\alpha}\ot\omega_\alpha)\Delta(h)
&=\frac1{c}\sum_{i,j}\lambda_i \omega_{\overline\alpha}(S(e_{ij}))\omega_\alpha(e_{ij})\\
&=\frac1{c}\sum_{i}\lambda_i=1.
\end{align*}
\ebew

Observe that this is a result already obtained in the original paper on discrete quantum groups, see Proposition 4.4 in  \cite{VD-discrete}.
\ssnl

We are now ready for the main property here. We refer to our papers on separability idempotents \cite{VD-si1,VD-si2}. 

\prop
The element $\Delta(h)$ is a separability idempotent in  $M(A\ot A)$.
\eprop

\bew
i) We first look at the setting in \cite{VD-si1}. For both the algebras $B$ and $C$ in that paper, we simply take $A$, and for $E$ we take $\Delta(h)$. The algebra $A$ is non-degenerate and $\Delta(h)$ is an idempotent in $M(A\ot A)$. Here we have $^*$-algebras and $\Delta(h)$ is self-adjoint. We know that $\Delta(h)(1\ot a)$ and $(a\ot 1)\Delta(h)$ both belong to $A\ot A$ and similarly on the other side. This shows that  Assumption 1.1 of \cite{VD-si1} are fulfilled. In particular $\Delta(h)$ is regular in the sense of that assumption.
\ssnl
ii) We claim that $\Delta(h)$ is full in the sense that the legs of $\Delta(h)$ are all of $A$. This is proven in the paper, but follows more easily from the equation $$H=\sum_i \lambda_i \,\Delta(h)(1_{\overline\alpha}\ot 1_\alpha)$$ Indeed, if $V$ is a subspace of $A$ satisfying $\Delta(h)(1\ot a)\subseteq V\ot A$ we get, with $a=e_{pq}$ from
\begin{equation*}
\sum_{ij} \lambda_i\, S(e_{ij})\ot e_{ji}e_{pq}=\sum_j \lambda_p S(e_{pj})\ot e_{jq}
\end{equation*}
that $S(e_{pj})\in V_{\overline\alpha}$ for all $j,p$. This implies that  $V_{\overline\alpha}=A_{\overline\alpha}$ for all $\alpha$ and hence $V=A$. By the way, from Proposition 1.3 in \cite{VD-si1} we can conclude that the legs of $\Delta(h)$ are all of $A$, now in the sense of our paper here. Also Proposition 1.4 in \cite{VD-si1} here says that $\Delta(h)(1\ot a)=0$ only if $a=0$, and similarly $(a\ot 1)\Delta(h)=0$ only if $a=0$. This we knew already because we have a multiplier Hopf algebra and $h\neq 0$.
\ssnl
iii) Finally we come to Definition 1.5 of \cite{VD-si1}. Since we know here already that 
\begin{equation*}
\Delta(h)(a\ot 1)=\Delta(h)(1\ot S(a)) \tussenen (1\ot a)\Delta(h)=(S(a)\ot 1)\Delta(h)
\end{equation*}
for all $a$, and because $S$ is bijective, we finally get that $\Delta(h)$ satisfies all the requirements and so it is a regular separability idempotent.
\ebew

The maps $S:B\to C$ and $S':C\to B$ we have in Proposition 1.7 of \cite{VD-si1} are precisely the antipode $S$ of $A$. If we look at the Proposition 1.8 in \cite{VD-si1}, we get here
\begin{align*}
m(S\ot \iota)(H(1\ot e_{pq}))
&=\sum_{i,j} \lambda_i S^2(e_{ij})e_{ji}e_{pq}\\
&=\sum_{j} \lambda_p S^2(e_{pj})e_{jq}\\
&=\sum_j \lambda_j e_{pj}e_{jq}=\sum_j \lambda_j e_{pq}
\end{align*}
and consequently $m(S\ot \iota)(\Delta(h)(1\ot a))=a$ for all $a$.
\ssnl
Similarly
\begin{align*}
m((\iota\ot S)((e_{pq}\ot\iota)H)
&=\sum_{i,j} \lambda_i\, e_{pq}S(e_{ij})S(e_{ji})\\
&=\sum_{i,j} \lambda_i\, e_{pq}S(e_{ji}e_{ij})\\
&=\sum_{i,j} \lambda_i\, e_{pq}S(e_{jj})\\
&=\sum_{i} \lambda_i\, e_{pq}
\end{align*}
and consequently also $m(\iota\ot S)\Delta(h))=a$ for all $a$.
\newpage
\bf The invariant integrals \rm
\nl
We consider the left and right integrals $\varphi$ and $\psi$, normalized such that $\varphi(h)=\psi(h)=1$. Then we have $\varphi\circ S=\psi$ as well as $\psi\circ S=\varphi$, see Proposition \ref{prop:1.12a}. 
\ssnl
We take a positive invertible element in $M(A)$ satisfying $aq=S^2(a)q$ for all $a$. Then we get the following formulas for the integrals. Recall that we use $\omega_\alpha$ for the standard trace on $A_\alpha$, see Notation \ref{notat:2.15a}.

\stel\label{stel:2.18a}
When $a\in A_\alpha$  we have
\begin{equation}
\varphi(a)=\omega_\alpha(q)\omega_\alpha(aq\inv).\label{eqn:2.4a}
\end{equation}
\estel

\bew
Define $\varphi$ on $A$ by $\varphi(a)=\omega_\alpha(q)\omega_\alpha(aq\inv)$ when $a\in A_\alpha$. We have to show that $(\iota\ot\varphi)\Delta(h)=1$. This means that we need $(\iota\ot\varphi)(\Delta(h)(1_{\overline\alpha}\ot 1))=1_{\overline\alpha}$ for all $\alpha$.
\ssnl
Now $\Delta(h)(1_{\overline\alpha}\ot 1)=\Delta(h)(1_{\overline\alpha}\ot 1_\alpha)$ and we can use the formula of Proposition \ref{prop:2.16a}. We get
\begin{align*}
(\iota\ot\varphi)(\Delta(h)(1_{\overline\alpha}\ot 1_\alpha))
&= \frac 1{c}\sum_{i,j} \lambda_i\, S(e_{ij})\varphi(e_{ji})\\
&= \frac 1{c}\sum_{i,j} \lambda_i\, S(e_{ij})c\omega_\alpha (e_{ji}q\inv)\\
&=\sum_i \lambda_i\, S(e_{ii})\lambda_i\inv\\
&=\sum_i S(e_{ii})=S(1_\alpha)=1_{\overline\alpha}.
\end{align*}
This proves the result. 
\ebew

Observe that the right hand side of Equation (\ref{eqn:2.4a}) does not depend on the choice of $q$. 

\stel\label{prop:2.19a}
For any $a\in A_\alpha$ we have
\begin{equation}
\psi(a)=\omega_\alpha(q\inv)\omega_\alpha(aq).\label{eqn:2.5a}
\end{equation}
\estel

\bew
i) Take $a\in A_{\overline\alpha}$. Then $S\inv(a)\in A_\alpha$ and we can apply the formula for $\varphi$. We find
\begin{align*}
\psi(a)=\varphi(S\inv(a))
&=\omega_\alpha(q)\omega_\alpha(S\inv(a)q\inv)\\
&=\omega_\alpha(q)\omega_\alpha(S\inv(S(q\inv)a))\\
&=\omega_\alpha(q)\omega_{\overline\alpha}(S(q\inv)a)\\
&=\omega_\alpha(q)\omega_{\overline\alpha}(aS(q\inv)).
\end{align*}
ii) Denote $q'=S(q\inv)$ so that $q=S\inv({q'}\inv)$ and $\omega_\alpha(q)=\omega_{\overline\alpha}({q'}\inv)$.
If we apply $S$ on the equation $aq=qS^2(a)$ we get $S(q)S(a)=S^2(S(a))S(q)$. This holds for all $a$ and because $S$ is bijective we also have $S(q)a=S^2(a)S(q)$ and $aS(q\inv)=S(q\inv)S^2(a)$. As $q'=S(q\inv)$ we get $aq'=q'S^2(a)$ for all $a$.
\ssnl
iii) We see that $\psi(a)=\omega_{\overline\alpha}({q'}\inv)\omega_{\overline\alpha}(aq')$ for $a\in A_{\overline\alpha}$. We can scale $q'$ to get a positive invertible element. Then the result follows for $\overline\alpha$ instead of $\alpha$.
\ebew

Also here, the right hand side of Equation (\ref{eqn:2.5a}) does not depend on the choice of $a$.
\nl
We now obtain formulas for the other data, the modular automorphism $\sigma$, the modular element $\delta$. We also verify some formulas relating these data.

\prop
The $\alpha$-component  $\delta_\alpha$ of the modular element $\delta$ is given by
\begin{equation}
\delta_\alpha=\frac{\omega_{\alpha}(q\inv)}{\omega_\alpha(q)}\,q_\alpha^2.\label{eqn:2.6b}
\end{equation}
\eprop

\bew
The results follows from the two previous ones because $\psi(a)=\varphi(a\delta)$ for all $a$ in $A$.
\ebew

The right hand side of  equality Equation (\ref{eqn:2.6b} does not depend on the choice of $q$ because is we replace $q_\alpha$ by a multiple of it, we get the same result.
\ssnl
Recall that $a\delta=\delta S^4(a)$ for all $a$, see Proposition \ref{prop:1.14c}. And indeed, because 
$aq_\alpha=q_\alpha S^2(a)$, we have also $aq_\alpha^2=q_\alpha^2 S^4(a)$ for $a\in A_\alpha$  so that the $\alpha$-component of $\delta$ has to be a scalar multiple of $q_\alpha^2$.
\ssnl

We can also verify the formulas
\begin{equation*}
(\varphi\ot\iota)\Delta(h)=\delta
\tussenen
(\iota\ot\psi)\Delta(h)=\delta\inv.
\end{equation*}
For the first one we have
\begin{align*}
(\varphi\ot\iota)(\Delta(h)(1\ot 1_\alpha))
&=\frac1{\omega_\alpha(q)}\sum_{i,j} \lambda_i\varphi(S(e_{ij}))e_{ji} \\
&=\frac1{\omega_\alpha(q)}\sum_{i,j}\lambda_i \psi(e_{ij})e_{ji} \\
&=\frac1{\omega_\alpha(q)}\sum_{i} \omega_\alpha(q\inv)\lambda_i^2e_{ii} \\
&=\frac{ \omega_\alpha(q\inv)}{\omega_\alpha(q)}\sum_{i}\lambda_i^2e_{ii}\\
&=\frac{ \omega_\alpha(q\inv)}{\omega_\alpha(q)}q_\alpha^2=\delta_\alpha.
\end{align*}

For the second equation we have
\begin{align*}
(\iota\ot\psi)(\Delta(h)(1_{\overline\alpha}\ot 1))
&=\frac1{\omega_\alpha(q)}\sum_{i,j} \lambda_i S(e_{(ij}) \psi(e_{ji} \\
&=\frac1{\omega_\alpha(q)}\sum_{i} \lambda_i S(e_{(ii})\omega_\alpha(q\inv)\lambda_i \\
&=\frac{\omega_\alpha(q\inv)}{\omega_\alpha(q)}S(q_\alpha^2)=S(\delta_\alpha)=\delta\inv_{\overline\alpha}.
\end{align*}

We immediately get the modular automorphisms from these formulas.

\ssnl
We have, for $a,b\in A_\alpha$, 
\begin{equation*}
\varphi(ab)=\omega_\alpha(q)\omega_\alpha(abq_\alpha\inv)=\omega_\alpha(q)\omega_\alpha(bq_\alpha\inv a)=\omega_\alpha(q)\varphi(bq_\alpha\inv a q_\alpha).
\end{equation*}
We see that $\sigma(a)=q_\alpha\inv aq_\alpha$ and this is in agreement with $\sigma(a)=S^2(a)$.
\ssnl
Similarly, from $\psi(a)=\omega_\alpha(q\inv)\omega_\alpha(aq)$ when $a\in A_\alpha$, we get
$\sigma'(a)=q_\alpha aq_\alpha\inv$. This is in agreement with $\sigma'(a)=S^{-2}(a).$
\nl
%
%

\section{\hspace{-17pt}. Duality between discrete and compact quantum groups} \label{s:dual} 

In this section, we will consider the dual of a discrete quantum group. We will follow the approach to duality as developed for general algebraic quantum groups in \cite{VD-alg}. The dual will be a compact quantum group and we know that all compact quantum groups arise in this way. By looking at a compact quantum group as the dual of a discrete quantum group, we can obtain results about compact quantum groups in a different and easier way. We refer to the appendix for more comments on this matter.
\nl
\bf The dual of a discrete quantum group \rm 
\nl
We start with a discrete quantum group as in Definition \ref{defin:1.2}. So we have a direct sum $A$ of full matrix algebras $(A_\alpha)_{\alpha\in I}$  and a coproduct $\Delta$ on $A$ as in Definition \ref{defin:1.1}. 
\ssnl
There are various ways to construct the dual, but in this case, the easiest way is by taking the \emph{reduced dual space} as follows.
\notat
Let $B_\alpha$ be the  dual $A'_\alpha$ of $A_\alpha$ and use the pairing notation. So we write $\langle a, b\rangle$ for the value of an element $b$ from $B_\alpha$ in the point $a$ of $A_\alpha$.  We use $B$ to denote the direct sum of the spaces $B_\alpha$. We have a pairing of $A$ with $B$ given by the formula $(a,b)\mapsto \sum_{\alpha\in I}\langle a_\alpha,b_\alpha\rangle$ where now $a=(a_\alpha)_\alpha$ and $b=(b_\alpha)_\alpha$. The sum is well-defined because we only have finitely many terms that are non-zero. We consider $B$ as a subspace of the dual space $A'$ of $A$.
\enotat

Remark that for $b$ in the direct product of the spaces $B_\alpha$, we still can define $\langle a,b\rangle$ by this formula when $a\in A$ because also then we only have finitely many terms in the sum that are non-zero. Any element of $A'$ is represented by an element in this direct product.
\ssnl

The space $B$ is made into an associative algebra by dualizing the coproduct of $A$. It is a $^*$-algebra. We make the following choices. We will motivate them later, see Remark \ref{opm:3.4}. 
\ssnl
In the following statement,  $z^-$  denotes the complex conjugate of the complex number $z$.

\prop\label{prop:3.1}
For $b,b'\in B $ we use  $\langle a, bb'\rangle=\langle\Delta(a),b'\ot b\rangle$ for $a\in A$ to define the product  $bb'\in A'$. For $b\in B$ we also define $b^*$  by $\langle a,b^*\rangle= \langle S(a^*),b\rangle^-$ for $a\in A$ These definitions make $B$ into a unital $^*$-algebra. 
\eprop
\bew
Because $b$ and $b'$ are reduced functionals, so that they are supported on only finitely many components of $A$, this product is well-defined in the dual space $A'$. It will again be in $B$. Associativity of this product is a consequence of coassociativity of the coproduct $\Delta$ on $A$. The counit $\varepsilon$ belongs to $B$ because it is a linear functional with support in the one-dimensional component $A_e$, spanned by the cointegral. It will be the unit in the algebra $B$.
\ssnl
The map $b\mapsto b^*$ is clearly involutive as $S(S(a^*)^*)=a$ for all $a\in A$. And because $\Delta$ is a $^*$-homomorphism and $\Delta\circ S=\zeta (S \ot S)\Delta$, we get $(bb')^*={b'}^*b^*$ for all $b,b'\in B$. As before, we used $\zeta$ for the flip map. Therefore $B$ is a $^*$-algebra.
\ebew

We also dualize the original product on $A$. The result is the following.

\prop\label{prop:3.2}
The $^*$-algebra $B$ is a Hopf $^*$-algebra for the coproduct $\Delta$ defined by
$\langle a\ot a',\Delta(b)\rangle=\langle aa',b\rangle$ for $a,a'\in A$ and $b\in B$.
\eprop

\bew
i) When $b$ has its support in one component, say $A_\alpha$, then $\Delta(b)$ will have support in $A_\alpha\ot A_\alpha$. Therefore $\Delta(b)\in B\ot B$. Coassociativity of $\Delta$ on $B$ follows from associativity of the product on $A$. We have $\Delta(1)=1\ot 1$ because $\varepsilon$ is a homomorphism on $A$. We get a homomorphism on $B$ because the coproduct $\Delta$ on $A$ is a homomorphism. Finally, because $S((aa')^*)=S(a'{}^*a^*)=S(a^*)S(a'{}^*)$, we get that this coproduct is also a $^*$-homomorphism. 
\ssnl
ii) The counit $\varepsilon$ on $B$ is obtained by evaluation in $1$. The antipode $S$ on $B$ is given by the formula $\langle a,S(b)\rangle=\langle S^{-1}(a),b\rangle$ for $a\in A$. 
\ebew

In the formula for the antipode on the dual, the inverse appears as a result of our conventions. We give more details in the following remark.
\opm\label{opm:3.4}
 i) In the theory of Hopf algebras, the conventions regarding duality are simple. The product on the dual is the dual of the coproduct and the coproduct on the dual is the dual of the product. This convention is completely symmetric. Unfortunately, in the operator algebra approach to quantum groups, things are a bit more complicated. Usually, the product on the dual is obtained by dualizing the coproduct (as for Hopf algebras) but the coproduct on the dual is obtained by flipping the coproduct coming from simply dualizing the product. This convention is of course not symmetric.
\ssnl
ii) As mentioned already, discrete quantum groups first appeared as duals of compact quantum groups (within the operator algebra context) and in the more recent literature, this point of view is still quite common. This implies that the coproduct on the discrete quantum group is flipped. We are now looking reversely at a compact quantum group as dual to a discrete quantum group. Therefore, and in order to be in accordance with the existing literature on this subject, we have chosen to define the product on the dual $B$ of $A$ by using the flipped coproduct and thus we arrive at the formula $\langle a,bb'\rangle=\langle\Delta(a),b'\ot b\rangle$ as in Proposition \ref{prop:3.1}. On the other hand, the coproduct is defined as dual to the product giving the formula $\langle a\ot a',\Delta(b)\rangle=\langle aa',b\rangle$ as in Proposition \ref{prop:3.2}.
\ssnl
iii) In this case, we also use $\Delta$ for the coproduct, $\varepsilon$ for the counit, $S$ for the antipode on $B$. 
\eopm

For convenience of the reader, below we collect the important formulas within this point of view. 

\opm\label{opm:3.5}
 For $a,a'\in A$ and $b,b'\in B$ we have
\begin{align} \langle a,bb'\rangle &=\langle \Delta(a),b'\ot b \rangle\label{eqn:3.1a} \\
         \langle aa',b \rangle&=\langle a\ot a',\Delta(b)\rangle \label{eqn:3.2a} \\
         \langle S(a),b\rangle&=\langle a, S^{-1}(b)\rangle\label{eqn:3.3a}\\
         \langle a, b^* \rangle&=\langle S(a^*),b\rangle^-\label{eqn:3.4a}\\
         \langle a^*,b\rangle&=\langle a,S(b)^*\rangle^-\label{eqn:3.5a}.
\end{align}
Equation (\ref{eqn:3.5a}) follows from the Equations (\ref{eqn:3.3a}) and (\ref{eqn:3.4a}). Indeed
 we obtain
\begin{align*}
\langle a^*,b\rangle
&=\langle a^*,b^{**}\rangle=\langle S(a),b^*\rangle^-\\
&=\langle a,S\inv (b^*)\rangle^-=\langle a, S(b)^*\rangle^-
\end{align*}
for all $a,b$.
\eopm
Again observe  the presence of the inverse of the antipode in Equation (\ref{eqn:3.3a}). This is one of the consequences of this convention. For the same reason, there is  the lack of symmetry between (\ref{eqn:3.4a}) and (\ref{eqn:3.5a}). We do not get the usual formulas for a pairing of Hopf $^*$-algebras or multiplier Hopf $^*$-algebras (as e.g.\ in \cite{Dr-VD}), but rather a modified form.

\nl
\bf The integral on the dual \rm
\ssnl
We now define the integral on the dual $(B,\Delta)$. Recall that $h$ is the cointegral in $A$.

\prop\label{prop:3.3}
The map $b\mapsto\langle h,b\rangle$ is a positive linear functional on $B$. It is both left and right invariant.
\eprop
\bew
Left and right invariance come from the fact that $ah=\varepsilon(a)h$ and $ha=\varepsilon(a)h$ for all $a$.
\ssnl
To prove positivity take any $b\in B$. There exists an element $a\in A$ such that $b=\varphi_A(a\,\cdot\,)$ where $\varphi_A$ is the left integral on $A$ satisfying $\varphi_A(h)=1$. This is possible because $\varphi_A$ is faithful and $b$ is supported in only finitely many components of $A$. Then using the Sweedler notation,
\begin{align*} 
\langle h,bb^*\rangle
&=\langle \Delta(h), b^*\ot b\rangle\\
&=\sum_{(h)} \langle h_{(1)},b^*\rangle \langle h_{(2)},b\rangle\\
&=\sum_{(h)} \langle h_{(1)},b^*\rangle \varphi_A(a h_{(2)})\\
&=\sum_{(h)} \langle S(a)h_{(1)},b^*\rangle \varphi_A(h_{(2)})\\
&=\langle S(a),b^*\rangle=\langle S(S(a)^*),b\rangle^-\\
&=\langle a^*,b\rangle^-=\varphi_A(aa^*)^-=\varphi_A(aa^*)
\end{align*}
and positivity follows from the positivity of $\varphi_A$. 
\ebew

So we see that the dual of the discrete quantum group (seen as a multiplier Hopf $^*$-algebra with positive integrals) is a Hopf $^*$-algebra with positive integrals. We also see from the above equation that the integral is faithful. This is in fact automatic (see e.g.\ \cite{VD-alg}). We will encounter a similar result with a similar argument in Proposition \ref{prop:3.13a} later in this section, as well in Proposition \ref{prop:4.6a} in the next section with a dual form.

\prop
The modular automorphism $\sigma$ of $\varphi_B$ is given by
\begin{equation*}
\langle a,\sigma(b)\rangle=\langle S^{-2}(a)\delta,b\rangle
\end{equation*}
where $a\in A$ and $b\in B$.
\eprop
\bew
i) Let $d=\varphi_A(a\,\cdot\,)$ with $a\in A$. Then for all $c\in B$,
\begin{align*}
\varphi_B(dc)=\langle h, dc\rangle
&=\sum_{(h)} \langle h_{(1)},c\rangle \langle h_{(2)},d\rangle\\
&=\sum_{(h)} \langle h_{(1)},c\rangle \varphi_A(ah_{(2)})\\
&=\sum_{(h)} \langle S(a)h_{(1)},c\rangle \varphi_A(h_{(2)})\\
&=\langle S(a),c\rangle.
\end{align*}
ii) On the other hand, for all $b\in B$, 
\begin{align*}
\varphi_B(bd)=\langle h, bd\rangle
&=\sum_{(h)} \langle h_{(1)},d\rangle \langle h_{(2)},b\rangle\\
&=\sum_{(h)}\varphi_A(a h_{(1)})\langle h_{(2)},b\rangle)\\
&=\sum_{(h)} \varphi_A(h_{(1)})\langle S^{-1}(a)h_{(2)},b\rangle)\\
&= \langle S^{-1}(a)\delta,b\rangle.
\end{align*}
iii) It follows that $\langle S(a),\sigma(b)\rangle=\langle S^{-1}(a)\delta,b\rangle$ for all $a$ and hence 
$$\langle a,\sigma(b)\rangle=\langle S^{-2}(a)\delta,b\rangle$$
 for all $a\in A$ and $b\in B$.
\ebew

Let us make a remark about this formula for the modular automorphism $\sigma$ of the integral $\varphi_B$ on $B$.

\opm
i) The formulas for our convention are obtained by taking the opposite product on the algebra $B$ while the coproduct remains the same. Consequently, the left and right integrals remain the same, as well as the modular element. On the other hand, the modular automorphisms should be replaced by their inverses.
\ssnl
ii) For the usual duality between algebraic quantum groups $A$ and $B$ we have the formula for $\delta_B$ given by 
$\langle a,\delta_B\rangle=\varepsilon(\sigma^{-1}(a))$. We use the extension of the pairing to $A\times M(B)$. Because in our case, the modular automorphism $\sigma$ is the same as the square of the antipode (see Proposition \ref{prop:2.10a}, we will have $\langle a,\delta_B\rangle =\varepsilon(a)$ and this means that $\delta_B=1$. This is in agreement with the fact that the left integral is also right invariant.
\ssnl
iii)  Again for the usual duality we have the formula
$\langle a,\sigma(b)\rangle=\langle S^2(a)\delta\inv,b\rangle$. Then for the inverse of the modular automorphism, we get
$\langle S^{-2}(a)\delta,\sigma\inv(b)\rangle$. This is in agreement with the formula we got in the previous proposition because, passing to our conventions, the antipode on $A$ remains the same, as well as the modular element, while the modular automorphism on $B$ has to be replaced by its inverse.
\eopm

Formulas of this type are found at various places in the literature. One possible reference is \cite{De-VD1} where they are obtained in the more general case of algebraic quantum hypergroups. They are also found in \cite{VD-part1}.
\ssnl
If we let $B$ act on the GNS-representation space associated with this integral, we can complete it to a C$^*$-algebra. Also the coproduct can be extended. One of the key features that makes this possible is the fact that the $^*$-algebra is represented by \emph{bounded} operators in the GNS-space associated with the positive integral. This is not true for general positive linear functionals on a $^*$-algebra. This result is already found in \cite{Ku-VD}, see also  \cite{DC-VD} and the more recent approach in \cite{VD-part2}.
\ssnl

In what follows we will mostly denote the integral on the dual compact quantum group $B$, that we obtained in Proposition \ref{prop:3.3}, by $\varphi_B$ and we will use $\varphi_A$ and $\psi_A$ for the left and the right integrals on the original discrete quantum group $A$. If no confusion is possible, we will drop these subscripts.
\nl
\bf Corepresentations of the dual \rm 
\nl
As one can expect, there is a one-to-one correspondence between corepresentations of the dual $(B,\Delta)$ and representations of the discrete quantum group $(A,\Delta)$. We make this precise and  examine it in detail.
\ssnl
Let $n$ be given and consider an element $u\in B\ot M_n(\mathbb C)$. Associate a linear map $\pi$ from  $A$ to $M_n(\mathbb C)$ by  $\pi(a)=(\omega_a\ot \iota)u$ where  $\omega_a=\langle a,\,\cdot\,\rangle$. In what follows, we use the leg-numbering notation as explained in the introduction.

\prop
The map $\pi$ is a homomorphism if and only if
$(\Delta\ot\iota)u=u_{13}u_{23}$. It is a non-degenerate homomorphism if and only if $u$ is invertible. In that case the inverse of $u$ is $(S\ot\iota)u$. It is a $^*$-representation if and only if $u$ is a unitary element in the $^*$-algebra $B\ot M_n(\mathbb C)$.
\eprop

\bew
i) Given $a,a'\in A$ we have
\begin{align*} \pi(aa')&=(\omega_{aa'}\ot \iota)u =((\omega_a\ot\omega_{a'}\ot\iota)((\Delta\ot\iota)u)\\
\pi(a)\pi(a')&=((\omega_a\ot\iota)u)((\omega_{a'}\ot\iota)u)=(\omega_a\ot\omega_{a'}\ot\iota)(u_{13}u_{23}).
\end{align*}
If $\pi(aa')=\pi(a)\pi(a')$ for all $a,a'$, it follows that $(\Delta\ot\iota)u=u_{13}u_{23}$. Clearly also the converse is true.
\ssnl
ii) Assume now that $\pi$ is a homomorphism such that $(\Delta\ot\iota)u=u_{13}u_{23}$ holds. If we apply the counit $\varepsilon$ of $B$ on the first factor of the equation $(\Delta\ot\iota)u=u_{13}u_{23}$, we see that $u=eu$ where $e=(\varepsilon\ot\iota)u$. Similarly if we apply it on the second factor. Then we get $ue=u$. 
If we apply $m(S\ot \iota)$ on the equation $(\Delta\ot\iota)u=u_{13}u_{23}$ we find
\begin{equation*}
((S\ot \iota)u)u=(\varepsilon\ot\iota)u=e
\end{equation*}
and similarly, if we apply $m(\iota\ot S)$ we get $(u(S\ot \iota)u)=e$. We see that $u$ is invertible if and only if $e=1$. This happens if and only if the representation $\pi$ is non-degenerate. 
\ssnl
iii) 
Finally, if $u$ is a unitary we find
\begin{align*} \pi(a)^* &=((\omega_a\ot\iota)u)^*=(\overline{\omega_a}\ot\iota)(u^*) \\
                  &=(\overline{\omega_a}\ot\iota)((S\ot 1)u).
\end{align*}
But we have 
$$\overline{\omega_a}(S(b))=\omega_a(S(b)^*)^-=\langle a,S(b)^*\rangle^-
     =\langle a^*,b\rangle.
$$
This means that $\overline{\omega_a}\circ S=\omega_{a^*}$
and so $\pi(a)^*=\pi(a^*)$. It is also easy to see that conversely, when $\pi(a)^*=\pi(a^*)$ for all $a$, we will  have that $u^*=(S\ot\iota)u$ so that $u$ is a unitary element in $B\ot M_n(\mathbb C)$.
\ebew

We also have a converse of the previous result.

\prop
For every $n$-dimensional representation $\pi$ of $A$, there is an element $u\in B\ot M_n(\mathbb C)$ such that $\pi(a)=(\omega_a\ot\iota)u$ for all $a\in A$.
\eprop
\bew
Let $\pi$ be an $n$-dimensional representation of $A$. Define linear functionals $u_{ij}$ on $A$ by $\langle a,u_{ij}\rangle=\pi(a)_{ij}$. We claim that these functionals must be elements in $B$. Then we put $u=\sum_{ij} u_{ij} \ot e_{ij}$ where $e_{ij}$ are matrix units. We clearly will have 
\begin{equation*}
(\omega(a\,\cdot\,)\ot \iota)u=\sum_{ij} \langle a,u_{ij}\rangle e_{ij}=\sum_{ij} \pi(a)_{ij} e_{ij}=\pi(a).
\end{equation*}
To prove the claim, observe that for only finitely many indices $\alpha\in I$ we can have $\pi(e_\alpha)\neq 0$. So the functionals $u_{ij}$ will be supported on finitely many components and hence they belong to $B$.
\ebew

This in the end shows that indeed, there is a one-to-one correspondence between $n$-dimensional representations of $A$ and elements in elements $u\in B\ot M_n(\mathbb C)$ satisfying $(\Delta\ot\iota)u=u_{13}u_{23}$.
\nl
Let us now have a look at the more general case of a \emph{possibly infinite-dimensional} corepresentation with coefficients in the C$^*$-algebra completion $\overline B$ of $B$ and not just in $B$. We now only look at unitary corepresentations. So in this case $u$ is in the multiplier algebra $M(\overline B\ot\Cal K(\Cal H))$ where $\Cal K(\Cal H)$ is the C$^*$-algebra of compact operators on a Hilbert space $\Cal H$. The C$^*$-tensor product is taken. It is still required that $(\Delta\ot\iota)u=u_{13}u_{23}$, this time in the multiplier algebra $M(\overline B\ot\overline B \ot\Cal K(\Cal H))$.
\ssnl
We can again associate a $^*$-representation of $A$ on the Hilbert space $\Cal H$ using the same formula, but we will need the following lemma to do this.

\lem
For all $a\in A$ the linear map $b\mapsto \langle a,b\rangle$ is continuous on $B$ and therefore can be extended to the C$^*$-algebra $\overline B$.
\elem
\bew
We have 
$$a=(\iota\ot\varphi_A)(\Delta(h)(a\ot 1))=(\iota\ot\varphi_A)(\Delta(h)(1\ot S(a))$$
so that $\langle a,b\rangle=\varphi_B(b'b)$ with $b'=\varphi_A(\,\cdot\,S(a))$. Because $b \mapsto \varphi_B(b'b)$ is continuous, this proves the result. 
\ebew

Then we can show the following.

\prop\label{prop:3.12a}
For a unitary corepresentation $u\in M(\overline B\ot\Cal K(\Cal H))$ we can define a non-degenerate $^*$-representation $\pi$ of $A$ on $\Cal H$ by $\pi(a)=(\omega_a\ot \iota)u$ where as before $\omega_a=\langle a,\,\cdot\,\rangle$, now extended to $\overline B$.
\eprop
\bew
Remark that $\pi(a)$ is defined as a multiplier of $\Cal K(\Cal H)$ and so $\pi(a)$ is defined in $\Cal B(\Cal H)$. The rest of the proof is essentially the same as for finite-dimensional corepresentations. 
\ebew

If we take any central projection $e$ of $A$ and if we cut down $\pi$ by $e$, we see that the resulting representation lives on a finite number of components of $A$. It follows that $u(1\ot \pi(e))$ is actually in the algebraic tensor product $B\ot \Cal B(\pi(e)\Cal H)$. Therefore, $u$ is a multiplier of the algebraic tensor product $B\ot \Cal K_0(\Cal H)$ where $\Cal K_0(\Cal H)$ denotes the $^*$-subalgebra of operators living on these finite-dimensional subspaces.  In particular, if $\Cal H$ is finite-dimensional, then $u\in B\ot \Cal B(\Cal H)$ and we get a corepresentation as we defined it earlier.
\nl
\bf The left regular corepresentation \rm 
\nl
Next we look at the left regular corepresentation of $B$ and show that it gives the GNS-representation of $A$ associated with the left integral $\varphi_A$. The unitary realizing the equivalence is the {\it Fourier transform}.
\ssnl
Consider the GNS-representation of $B$ associated with the integral and denote by $\Lambda_B: B \to \Cal H_B$ the associated canonical map. We let $B$ act directly on the Hibert space $\Cal H_B$ so that $\Lambda_B(bb')=b\Lambda_B(b')$ for all $b,b'\in B$. Similarly, take the GNS-representation of $A$ for the left integral $\varphi_A$ on $A$ and denote the canonical map by $\Lambda_A: A \to \Cal H_A$.

\prop\label{prop:3.13a}
With the above notations, there is a unitary operator $F:\Cal H_A \to \Cal H_B$ given by $F\Lambda_A(a)=\Lambda_B(b)$ were $a\in A$ and $b=\varphi_A(S^{-1}(\,\cdot\,)a)$. If $\pi$ is the $^*$-representation of $A$ on $\Cal H_B$ determined by the left regular corepresentation of $B$, then $F^*\pi(\,\cdot\,)F$ is the GNS representation of $A$ on $\Cal H_A$.
\eprop
\bew
i) First, we will show that $\varphi_B(b^*b)=\varphi_A(a^*a)$ if $a\in A$ and $b=\varphi_A(S^{-1}(\,\cdot\,)a)$. This is a modified form of the formula found in the proof of Proposition \ref{prop:3.3}. The argument is essentially the same. We have
\begin{align*} 
\varphi_B(b^*b) 
&= \langle h, b^*b\rangle =\langle \Delta(h),b\ot b^*\rangle \\
&=\sum_{(h)} \langle h_{(1)},b\rangle \langle h_{(2)},b^* \rangle \\
 &=\sum_{(h)} \varphi_A(S\inv(h_{(1)})a) \langle h_{(2)},b^* \rangle \\
 &=\sum_{(h)} \psi_A (S(a)h_{(1)}) \langle h_{(2)},b^* \rangle \\
 &=\sum_{(h)} \psi_A (h_{(1)}) \langle ah_{(2)},b^* \rangle \\
 &=\langle a,b^*\rangle=\langle S(a^*),b\rangle^-\\
  &=\varphi_A(a^*a)^-=\varphi_A(a^*a).
\end{align*}
This implies the existence of the unitary operator $F:\Cal H_A \to \Cal H_B$ given by $F\Lambda_A(a)=\Lambda_B(b)$.
\ssnl
ii) The left regular corepresentation $u$ is defined by 
$$u^*(\eta\ot\Lambda_B(b))=\sum_{(b)}b_{(1)}\eta\ot\Lambda_B(b_{(2)}),$$
see the papers on locally compact quantum groups, mentioned in the introduction. Because $u^*=(S\ot \iota)u$ we get
$$u(\eta\ot\Lambda_B(b))=\sum_{(b)}S^{-1}(b_{(1)})\eta\ot\Lambda_B(b_{(2)}).$$
Therefore, with $\pi$ as in Proposition \ref{prop:3.12a} we find
$$\pi(a')\Lambda_B(b)=\sum_{(b)}\langle a',S^{-1}(b_{(1)}) \rangle\Lambda_B(b_{(2)})$$
when $a'$ is another element in $A$. Now
\begin{align*} \sum_{(b)}\langle a',S^{-1}(b_{(1)}) \rangle &\langle x,b_{(2)}\rangle
      =\langle S(a')x,b\rangle \\
      &= \varphi_A(S^{-1}(S(a')x)a) =\varphi_A(S^{-1}(x)a'a)
\end{align*}
and therefore
$$\sum_{(b)}\langle a',S^{-1}(b_{(1)}) \rangle \Lambda_B(b_{(2)})
      =F\Lambda_A(a'a)=F a'\Lambda_A(a).$$
Because $\Lambda_B(b)=F \Lambda_A(a)$, we get indeed 
$F^*\pi(a')F=a'$ on $\Cal H_A$.
\ebew

\snl
\bf More formulas \rm 
\nl
Finally, we look at some  important formulas we find in the theory of compact quantum groups and see what the corresponding results are for the discrete quantum groups.
\ssnl
The orthogonality relations of irreducible unitary corepresentations yield an operator $Q$ in $M(A)$ by the formula
$$\varphi_B(u^\alpha_{ip}(u^\beta_{jq})^*)=\delta_{\alpha\beta}\delta_{ij}Q^\alpha_{qp}.$$
In the following formula, we find that this element is nothing else but (a scalar multiple of) $\delta^{-\frac12}$ where now $\delta$ is the modular element of the discrete quantum group, relating the left and the right integrals.

\prop
Choose an index $\alpha$ and denote by $v$ the associated unitary corepresentation of $B$. The
\begin{equation*}
\varphi_B(v_{ij}(v_{ik})^*)=(q_\alpha)_{jk}.
\end{equation*}
\eprop

\bew
We know that $\varphi_A(aa^*)=\varphi_B(bb^*)$ if $b=\varphi_A(a\,\cdot\,)$ (see Proposition \ref{prop:3.3}). So we have to find the elements $a_{ij}\in A_\alpha$ such that $v_{ij}=\varphi_A(a_{ij}\,\cdot\,)$. This means that we must have
$$\varphi_A(a_{ij}x)=\langle x,v_{ij}\rangle = x_{ij}$$
when $x\in A_\alpha$.
So 
$$x_{ij}
=\omega_\alpha(a_{ij}xq_\alpha\inv)
=\omega_\alpha(q_\alpha\inv a_{ij}x).$$
Therefore $q_\alpha\inv a_{ij}=e_{ji}$ where we have matrix units in $A_\alpha$.
Then 
\begin{align*} \varphi_B(v_{ij}v_{ik}^*)
 & =\varphi_A(a_{ij}a_{ik}^*)\\
 &=\omega_\alpha(q_\alpha e_{ji} e_{ik} q_\alpha q_\alpha\inv)\\    
 &=\omega_\alpha(q_\alpha e_{jk})
 =(q_\alpha) _j \delta_{jk}.
\end{align*}

\vskip -15pt
\ebew

%
%

\section{\hspace{-17pt}. Discrete and compact groups} \label{s:examples} 

In this section, we will discuss  two special cases to illustrate the results of the previous sections. We begin with the motivating case of a discrete quantum group arising from a discrete group. Most of this is trivial but we include it for completeness. There are also some consequences of the choice we made as in Remark \ref{opm:3.5} in the previous section. 
\ssnl
We also look at the dual of a compact group.  In \cite{VD-suq2} we treat a more advanced example in greater detail.
\nl
\bf The discrete quantum group arising from a discrete group \rm
\nl

Recall Example \ref{voorb:2.3a} from Section \ref{s:discr}. In that example, $G$ is any group. A discrete quantum space (as in the beginning of that section) is associated. The index set  is identified with the set $G$ and let $n(p)=1$ for all $p$ in $G$. So every component $A_p$ is equal to $\mathbb C$. The direct sum $A$ of these components is nothing else but the algebra $K(G)$ of all complex functions on $G$ with finite support and  pointwise operations. A coproduct on $A$ is defined by $\Delta(f)(p,q)=f(p,q)$.  Then the pair $(A,\Delta)$ is a discrete quantum group as in Definition \ref{defin:1.2}.
\ssnl 
Here the algebra $A$ is abelian. We also have the converse result.

\prop 
If $(A,\Delta)$ is a discrete quantum group and if $A$ is abelian, then it is the discrete quantum group associated to a group as above.
\eprop

\bew
Now identify the set $G$ with the index set $I$ of the discrete quantum group $(A,\Delta)$. For all pairs $p,q\in G$ we consider the homomorphism $f\mapsto \Delta(f)(p,q)$ form $A$ to $\mathbb C$. It is given by evaluation in a point of $G$ that we denote by $pq$. This gives an associative product on $G$ because $\Delta$ is coassociative. There is an identity $e$ for the  multiplication  defined by $f(e)=\varepsilon(f)$. And every element $p\in G$ has an inverse $p\inv$ given by $f(p\inv)=S(f)(p)$ for all $f\in K(G)$. 
\ebew

In this case the cointegral $h$ is the function $\delta_e$  with value $1$ in $e$ and $0$ everywhere else. For $\Delta(h)$ we can write
\begin{equation*}
\Delta(h)=\sum_p \delta_p\ot \delta_{p\inv},
\end{equation*}
where just like for $\delta_e$, the function $\delta_p$ takes the value $1$ in the point $p$ and $0$ everywhere else. 
\ssnl

The left  integral $\varphi_A$ is given by $\varphi(f)=\sum_p f(p)$. It is also a right integral. 
 It is straightforward to verify left and right invariance of $\varphi_A$. Indeed we have 
\begin{align*}
(\iota\ot\varphi_A)\Delta(h)&=\sum_p \delta_p \varphi_A(\delta_{p\inv})=\sum_p\delta_p=1,\\
(\varphi_A\ot \iota)\Delta(h)&=\sum_p \varphi_A(\delta_p) \delta_{p\inv}=\sum_p  \delta_{p\inv}=1.
\end{align*}

These equations hold in the multiplier algebra of $K(G)$.
\ssnl
As the algebra is abelian, these functionals are traces and the modular automorphisms are trivial.  Also the modular element $\delta$ is trivial because the left and right integrals coincide. The square $S^2$ of the antipode is the identity map. 
\ssnl
Let us have a quick look at the formula in Proposition \ref{prop:2.11a}. For the standard trace on the components, we simply have $\omega(\delta_p)=1$ for all $p$.  For the element $q$, defined in Proposition \ref{prop:2.11a}, we get here 
\begin{equation*}
q=(\omega\ot\iota)\Delta(h)=\sum_p \omega(\delta_p) \delta_{p\inv}=\sum_p \delta_{p\inv}=1.
\end{equation*}
 Because we took the normalized trace, we get $q=1$ in  this case as expected. 
\nl
We now pass to the results of Section \ref{s:dual} for this example. 
\ssnl
As before, let $G$ be a group and $(K(G),\Delta)$   the associated discrete quantum group. We also consider the Hopf $^*$-algebra $(\mathbb C G,\Delta)$. Here $\mathbb C G$ is the group algebra and $\Delta$ the coproduct on $\mathbb C G$ given by $\Delta(\lambda_p)=\lambda_p\ot\lambda_p$. We use $p\mapsto \lambda_p$ for the canonical embedding of the group in the group algebra. 
\prop
Define a linear map $\gamma$ from $\mathbb C G$ to the dual of $K(G)$ by 
$$\gamma(\lambda_p)=\varphi_A(\,\cdot\,\delta_{p\inv}).$$
 Then $\gamma $ is an isomorphism from $(\mathbb C G,\Delta)$ to the  dual  of $(K(G),\Delta)$. 
\eprop
\bew
It is clear that the map $\gamma$ is a bijective map from the space $\mathbb C G$ to the underlying space of the dual of the discrete quantum group $(K(G),\Delta)$. 
\ssnl
Consider the pairing of $K(G)$ with its dual. Then we have
\begin{equation*}
\langle f, \gamma(\lambda_p)\rangle=\varphi(f\delta_{p\inv})=f(p\inv)
\end{equation*}
for all $f\in K(G)$ and $p\in G$.
\ssnl
i) Given $f\in K(G)$ and $p,q\in G$ we have
\begin{align*}
\langle f,\gamma(\lambda_{pq})\rangle
&=f((pq)\inv)=f(q\inv p\inv)\\
&=\Delta(f)(q\inv,p\inv)\\
&=\langle \Delta(f),\gamma(\lambda_q)\ot \gamma(\lambda_p) \rangle\\
&=\langle f, \gamma(\lambda_p)\gamma(\lambda_q) \rangle\
\end{align*}
Remember the convention about the product on the dual algebra (cf. Equation (\ref{eqn:3.1a}) in the previous section). So we get $\gamma(\lambda_p\lambda_q)=\gamma(\lambda_p)\gamma(\lambda_q)$ and $\gamma$ is a homomorphism from the group algebra $\mathbb C G$ to the dual of $K(G)$. 
\ssnl
ii) Let $f\in K(G)$ and $p\in G$. Then 
\begin{align*}
\langle f,\gamma(\lambda_p^*)\rangle
&=\langle f, \gamma(\lambda_{p\inv})\rangle=f(p),\\
\langle f,\gamma(\lambda_p)^* \rangle
&=\langle S(f^*),\gamma(\lambda_p)\rangle^-\\
&=(S(\overline f)(p\inv)^-=f(p)
\end{align*}
and we see that $\gamma(\lambda_p^*)=\gamma(\lambda_p)^*$. Here we have used Equation (\ref{eqn:3.4a}) from the previous section. Therefore $\gamma$ is a $^*$-isomorphism from $\mathbb C G$ to the dual of $K(G)$.
\ssnl
iii) Given $f,g\in K(G)$ and $p\in G$ we have (using Equation (\ref{eqn:3.2a}),
\begin{align*}
\langle fg,\gamma(\lambda_p)\rangle
&=(fg)(p\inv)=f(p\inv)g(p\inv)\\
&=\langle f\ot g,\gamma(\lambda_p)\ot \gamma(\lambda_p)\rangle\\
&=\langle f\ot g,(\gamma\ot \gamma)\Delta(\lambda_p)\rangle.
\end{align*}
This means that $\Delta(\gamma(\lambda_p))=(\gamma\ot \gamma)\Delta(\lambda_p)$. This completes the proof.
\ebew

The above result gives the pairing of $K(G)$ and $\mathbb C G$ as in the following remark.

\opm Consider the $^*$-isomorphism $\gamma$ from $\mathbb C G$ to $\widehat{K(G)}$ given in the previous proposition. It induces a pairing of $K(G)$ with $\mathbb C G$ by the formula
\begin{equation*}
\langle f,\lambda_p\rangle=\langle f,\gamma(\lambda_p)\rangle=f(p\inv).
\end{equation*}
Again because of our conventions, this is not the usual pairing of these two algebras.
\eopm

Let us also consider the dual integral. The right integral $\psi_B$ on the dual is defined by $\psi_B(b)=\varepsilon(f)$ when $b=\varphi_A(\,\cdot\, f)$. Hence we get 
\begin{equation*}
\psi_B(\gamma(\lambda_p))=\varepsilon( \delta_{p\inv})=\delta(e,p\inv)=\delta(e,p)
\end{equation*}
and indeed $\psi_B\circ \gamma$ is the integral on the group algebra $\mathbb C G$.
\nl
\bf The dual of the group algebra $\mathbb C G$ \rm
\nl
In this item, we start with a group $G$ and its group algebra $\mathbb C G$. We also use $B$ for this algebra. We consider the coproduct $\Delta$ defined by $\Delta(\lambda_p)=\lambda_p\ot \lambda_p$. Recall that $p\mapsto \lambda_p$ is the canonical embedding of $G$ in $\mathbb C G$. This gives a Hopf $^*$-algebra with integrals. The left integral $\varphi_B$ is also a right integral and $\varphi_B(\lambda_p)=0$ except for $p=e$. Recall that we use $e$ for the identity of the group $G$. If we normalize the left integral $\varphi_B$  by letting $\varphi_B(\lambda_e)=1$, we get a positive integral.
\ssnl
Indeed, we have for $b=\sum_p b(p)\lambda_p$ that
\begin{equation*}
\varphi_B(b^*b)=\sum_{p,q} \overline{b(p)}b(q) \varphi(\lambda_{p\inv q})=\sum_{p} \overline{b(p)}b(p). 
\end{equation*}
\ssnl
Then we can apply the general duality theory of algebraic quantum groups and consider the dual $(\widehat B,\widehat\Delta)$ of this Hopf $^*$-algebra. Elements in $\widehat B$ are linear functionals on $A$ of the form $x\mapsto \varphi(xa)$ where  $a\in A$. The product on $\widehat B$ is dual to the coproduct on $A$ while, with our convention here, the coproduct on $\widehat B$ is now twisted. 
\ssnl
In what follows, we will use $A$ for the algebra $\widehat B$ and consider the pairing notation $(a,b)\mapsto \langle a,b\rangle$ where $a\in A$ and $b\in B$. With these conventions, the previous statement gives, more precisely, the following result.

\prop 
Define $\rho:K(G)\to A$ by 
$\rho(\delta_p)=\varphi_B(\,\cdot\,\lambda_p)$ for all $p\in G$. Then $\rho$ is an isomorphism of the discrete quantum group $(K(G),\Delta)$ to the dual  of the  Hopf $^*$-algebra $(\mathbb C G,\Delta)$.
\eprop
\bew

It is clear that $\rho$ is an isomorphism of vector spaces. For any $f\in K(G)$ and $p\in G$ we have
\ssnl
\begin{equation*}
\langle \rho(f),\lambda_p\rangle
=\sum_q f(q)\varphi_B(\lambda_p\lambda_q)
=\sum_q f(q)\varphi_B(\lambda_{pq})=f(p\inv).
\end{equation*}
\ssnl
ii) For any pair of functions $f,g\in K(G)$ we have, for all $p\in G$ that
\begin{align*}
\langle \rho(fg),\lambda_p\rangle
&=(fg)(p\inv)=f(p\inv)g(p\inv)\\
&=\langle \rho(f),\lambda_p\rangle\langle \rho(g),\lambda_p \rangle\\
&=\langle \rho(f)\ot\rho(g),\lambda_p\ot\lambda_p\rangle\\
&=\langle \rho(f)\ot\rho(g),\Delta(\lambda_p)\rangle\\
&=\langle \rho(f)\rho(g),\lambda_p\rangle
\end{align*}
and so $\rho(fg)=\rho(f)\rho(g)$ so that $\rho$ is a homomorphism. We have used Equation (\ref{eqn:3.2a}).
\ssnl
ii) We also have, using Equation (\ref{eqn:3.5a}),
\begin{align*}
\langle \rho(f^*), \lambda_p\rangle&=f^*(p\inv)=\overline f(p\inv)\\
\langle \rho(f)^*,\lambda_p\rangle&=\langle f,S(\lambda_p)^*\rangle^-\\
&=\langle f,\lambda_p\rangle^-=\overline f(p\inv)
\end{align*}
and we see that 
because $\rho(f^*)=\rho(f)^*$ . So we have that $\gamma$ is a $^*$-homomorphism
\ssnl
iii) Finally we show that $\Delta(\rho(f))=(\rho\ot\rho)\Delta(f)$ for all $f\in K(G)$. To show this, take $p,q\in G$. On the one hand we have, using Equation (\ref{eqn:3.1a}),
\begin{equation*}
\langle \Delta(\rho(f)),\lambda_p\ot\lambda_q \rangle
=\langle \rho(f),\lambda_q\lambda_p\rangle=\langle \rho(f),\lambda_{qp} \rangle= f(p\inv q\inv) 
\end{equation*}
while on the other hand
\begin{equation*}
\langle (\rho\ot\rho)\Delta(f), \lambda_p\ot\lambda_q\rangle
=\Delta(f)(p\inv,q\inv)=f(p\inv q\inv).
\end{equation*}
This completes the proof.
\ebew

This $^*$-isomorphism induces the same pairing between $K(G)$ and $\mathbb C G$ we have in the previous item. Indeed 
\begin{equation*}
\langle f,\lambda_p\rangle=\langle \rho(f),\lambda_p\rangle=f(p\inv).
\end{equation*}
In this equation, the left hand side is the pairing of $K(G)$ and $\mathbb C G$ as it induced by this isomorphism from the pairing of the dual of $\mathbb C G$ with $\mathbb C G$ in the middle expression. The last equality is by the definition of $\rho$ as we have seen in the beginning of the previous proof.

\newpage 
\bf The dual of a compact quantum group \rm
\nl
Next we start with a more general setting than in the previous item. We now consider any Hopf $^*$-algebra with positive integrals. These are the so-called $^*$-algebraic quantum groups of compact type, or shortly, as we will call them here, the compact quantum groups. 
\ssnl
So in what follows $(B,\Delta)$ is a Hopf $^*$-algebra with a positive integral $\varphi_B$. From the positivity of the integral, we must have that $\varphi_B(1)\geq 0$ and we can assume that $\varphi_B(1)=1$. It also follows that $\varphi_B$ is not only left invariant, but also right invariant.
\ssnl
We will now show that the dual $(A,\Delta)$ of $(B,\Delta)$, in the sense of duality of algebraic quantum groups, is a discrete quantum group as defined in Definition \ref{defin:1.2} of Section \ref{s:defin}. We will however work within the conventions used in Section \ref{s:dual}, cf.\ the formulas in Remark \ref{opm:3.5}.
\ssnl
 First we have the existence of a cointegral.

\prop\label{prop:4.5a}
The integral $\varphi_B$ on $B$ is a left cointegral in $A$. It is a self-adjoint idempotent.
\eprop
\bew
We define $h$ in $A$ by $\varphi_B$. This belongs to the dual of $B$ because $B$ has an identity.
\ssnl
i) Let $a\in A$.
Then for all $y\in B$, using that $\varphi_B$ is right invariant,  we have
\begin{align*}
\langle ah, y\rangle 
&=\langle a\ot h, \Delta(y),\rangle\\
&=\langle a,(\iota\ot \varphi_B)\Delta(y)\rangle\\
&=\varphi_B(y)\langle a,1\rangle=\langle h,y\rangle\, \varepsilon(a)
\end{align*}
and we see that $ah=\varepsilon(a)h$. Similarly, because $\varphi_B$ is also right invariant, we get $ah=\varepsilon(a)h$. 
\ssnl
ii) For the adjoint $h^*$ of $h$ in $B$ we get
\begin{equation*}
\langle h^*,y\rangle=\langle h,S(x)^*\rangle^-=\varphi_B(S(y)^*))^-=\varphi_B(S(y))=\varphi_B(y)
\end{equation*}
and $h^*=h$. We have used that $\varphi_B\circ S=\varphi_B$ here and that $\varphi_B$ is self-adjoint.
\ssnl
iii) For $h^2$  we get $\varepsilon(h)h$ and 
$\varepsilon(h)=\langle 1,h\rangle=\varphi_B(1)=1.$
\ebew

To conclude that the dual $(A,\Delta)$ of $(A,\Delta)$ is a discrete quantum group, we  use Theorem \ref{stel:2.8a} that we already proved in Section \ref{s:defin}.  For this, we need that the dual $A$ is an operator algebra. This will follow from the positivity of the integral. For this we can use the general result that the dual has a positive integral if the original $(B,\Delta)$ has a positive integral. 
\ssnl
We include the result  below and we keep using our conventions from Section \ref{s:dual}. For these conventions, as the coproduct on the dual $A$ is flipped  (cf.\ Equation (\ref{eqn:3.1a}), now the left integral $\varphi_A$ is  given by the formula $\varphi_A(a)=\varepsilon(b)$ when $a=\varphi_B(\,\cdot\,b)$ with $b\in A$.

\prop\label{prop:4.6a}
For any $a\in A$ of the form $\varphi_B(\,\cdot\, b)$ we get $\varphi_A(a^*a)=\varphi_B(b^*b)$.
\eprop

\bew
For all $y$ in $B$ we have
\begin{align*}
\langle a^*a,y\rangle
&=\langle a^*\ot a,\Delta(y)\rangle\\
&=\sum_{(y)} \langle a^*,y_{(1)}\rangle \langle a,y_{(2)} \rangle \\
&=\sum_{(y)} \langle a^*,y_{(1)}\rangle\varphi_B(y_{(2)}b) \\
&=\sum_{(b)} \langle a^*,S\inv(b_{(1)})\rangle \varphi_B(yb_{(2)})
\end{align*}
and it  follows that
\begin{equation*}
a^*a=\sum_{(b)} \langle a^*,S\inv(b_{(1)})\rangle \varphi_B(\,\cdot\,b_{(2)}).
\end{equation*}
Then
\begin{align*}
\varphi_A(a^*a)
&=\sum_{(b)} \langle a^*,S\inv(b_{(1)})\rangle \varepsilon(b_{(2)})\\
&=\langle a^*,S\inv(b) \rangle=\langle a,b^*\rangle^-\\
&=\varphi_B(b^*b)^-=\varphi_B(b^*b).
\end{align*}
\ebew
We have use the Equations (\ref{eqn:3.2a}) and (\ref{eqn:3.5a}) from Section \ref{s:dual}, as well as the formula
\begin{equation*}
S((\iota\ot\varphi_B)(\Delta(y)(1\ot b)))=(\iota\ot\varphi_B)((1\ot y)\Delta(b)).
\end{equation*}
This formula is true for any left integral on an algebraic quantum group. It is the dual form of the equality 
\begin{equation*}
(1\ot a)\Delta(h)=(S(a)\ot 1)\Delta(h)
\end{equation*}
as we explain further in Remark \ref{opm:4.8}.

\ssnl
Compare the above argument with the proof of Proposition \ref{prop:3.3}. There we get the positivity of the integral on $B$ from the positivity of the integral on $A$. Both results are essentially based on the same formula. This formula is like the Plancherel  formula for locally compact abelian groups, see more details about this in \cite{VD-part2}.

\ssnl
Now we can complete the proof of the following result.

\stel
The dual of a Hopf $^*$-algebra with a positive integral is a discrete quantum group.
\estel

\bew
Take a Hopf $^*$-algebra $B$ with a positive integral $\varphi_B$ as before. Let $A$ be the dual multiplier Hopf $^*$-algebra in the sense of duality of algebraic quantum groups. From Proposition \ref{prop:4.5a} we know that $A$ has a cointegral that is a self-adjoint idempotent.  From Proposition \ref{prop:4.6a} we know that $A$ again has a positive integral. The $A$ is an operator algebra and from from Theorem \ref{stel:2.8a} we get that $A$ is a direct sum of finite-dimensional matrix algebras. Then it is a discrete quantum group by Definition \ref{defin:1.2}
\ebew

We are using that a multiplier Hopf $^*$-algebra with a positive integral is an operator algebra. In this case, because of the structure of $A$, this is an easy result. In general, it uses that the algebra is represented by bounded operator in the GNS-representation of $A$ induced by the positive functional. See e.g.\ \cite{VD-part2} for details.

\opm\label{opm:4.8}
Take $x,y$ in $B$ and $a\in A$. We then have on the one hand
\begin{align*}
\langle (1\ot a)\Delta(h), x\ot y \rangle
&=\sum_{(h)} \langle h_{(1)},x\rangle \langle a h_{(2)},y\rangle\\
&=\sum_{(h),(y)} \langle h_{(1)},x\rangle \langle a,y_{(1)}\rangle \langle h_{(2)},y_{(2)}\rangle\\
&=\sum_{(y)} \langle h,y_{(2)}x\rangle \langle a,y_{(1)}\rangle \\
&=\langle a\ot h ,\Delta(y)(1\ot x) \rangle
\end{align*}
This is the pairing of $a$ with $(\iota\ot\varphi_B)(\Delta(y)(1\ot x))$. On the other hand
\begin{align*}
\langle (S(a)\ot 1)\Delta(h), x\ot y \rangle
&=\sum_{(h)} \langle S(a)h_{(1)},x\rangle \langle h_{(2)},y\rangle\\
&=\sum_{(h),(x)} \langle S(a),x_{(1)}\rangle \langle h_{(1)},x_{(2)}\rangle \langle h_{(2)},y\rangle\\
&=\sum_{(x)} \langle S(a),x_{(1)}\rangle \langle h,yx_{(2)}\rangle \\
&=\langle S(a)\ot h ,(1\ot y)\Delta(x) \rangle.
\end{align*}
This is the pairing of $S(a)$ with $(\iota\ot\varphi_B)((1\ot y)\Delta(x))$. We see that the equality 
$$(1\ot a)\Delta(h)=(S(a)\ot 1)\Delta(h)$$
 in $A\ot A$ is nothing else but the adjoint of the formula
\begin{equation*}
(\iota\ot\varphi_B)(\Delta(y)(1\ot x))=S\inv((\iota\ot\varphi_B)((1\ot y)\Delta(x)))
\end{equation*}
in  $B$. We have used the formulas we have in Remark \ref{opm:3.5}. 
\eopm

%
%
%

\section{\hspace{-17pt}. Conclusion and further research}\label{s:concl} 

In this paper, we have given a treatment of discrete quantum groups. It is different from the original one in \cite{VD-discrete}. The study of discrete quantum groups started before other more general developments in the theory of locally compact groups. We felt an update was needed with references to these newer developments. 
\ssnl
An important role is played by the element $\Delta(h)$. It is a separability idempotent in the sense of \cite{VD-si1} and many of its properties are in fact coming from this fact. 
\ssnl
We consider the dual of a discrete quantum group in the sense of duality of algebraic quantum groups as in \cite{VD-alg}. It is a compact quantum group and properties of the compact quantum group follow within this point of view of properties of discrete quantum groups.
\ssnl
In another paper on the subject \cite{VD-suq2}, the discrete quantum group of the quantized enveloping algebra of the Lie algebra of $SU(2)$ is studied in this setting. It illustrates many of the general results as described in this paper. It could be a challenging project to study more of such concrete examples in this framework.
\ssnl
A review of the more recent results on discrete quantum groups, obtained by various authors but with the approach of considering discrete quantum groups as duals of compact quantum groups also seems a desirable and interesting project.

%
%

\section*{\hspace{-17pt} Appendix. Different approaches to discrete quantum groups}\label{s:appA}

In this appendix, we will discuss and compare different approaches to discrete quantum groups and we will place our treatment in this context. We will also mention differences in conventions used within these different approaches.
\nl
\bf Discrete quantum groups as duals of compact quantum groups \rm
\nl
Compact quantum groups in the operator algebra framework are introduced and studied by Woronowicz. First he only considered a restricted class of such compact quantum groups in \cite{Wo2} while later, he studied them from a more general point of view \cite{Wo3}. This last paper only appeared in 1995, but it was available as a preprint many years before that. In this paper, we have used that a compact quantum group is in fact just a Hopf $^*$-algebra with a positive integral.
\ssnl
It seems appropriate to start this discussion with the work of Podle\'s and Woronowicz as they were the first to consider discrete quantum groups in an operator algebraic setting. We look at the Sections 2 and 3 in their paper on the quantum deformation of the Lorentz group \cite{Po-Wo}. 
\ssnl
In Section 2, they first review the theory of compact matrix quantum groups as developed in \cite{Wo2}. The underlying C$^*$-algebra is denoted by $A_c$ and we will denote the coproduct by $\Delta_c$. Then they use the results about corepresentations of $(A_c,\Delta_c)$ to construct the underlying algebra $A_d$ of the dual discrete quantum group. It is defined as a C$^*$-direct sum of matrix algebras. Corepresentation theory of $A_c$ is also used to construct the 'left regular corepresentation' $u$ in the multiplier algebra $M(A_d\ot A_c)$ of the C$^*$-tensor product of $A_d$ with $A_c$. 
\ssnl
The first main result in \cite{Po-Wo} is Theorem 3.1. The coproduct on the dual $A_d$ is constructed, together with a counit and an antipode satisfying the standard equalities. The coproduct is defined as a non-degenerate $^*$-homomorphism $\Delta_d:A_d\to M(A_d\ot A_d)$ by the formula $(\Delta_d\ot\iota)u=u_{23}u_{13}$. This means that they work with the 'flipped' coproduct on the dual, i.e.\ the usual convention in the operator algebra approach. The counit is a $^*$-homomorphism from $A_d$ to $\Bbb C$ (and so everywhere defined and bounded). The antipode however is not defined everywhere but only on the algebraic direct sum, denoted by $A_{d0}$, of the components of $A_d$ and maps this dense $^*$-subalgebra into itself.
\snl
In their Theorem 3.3, they obtain the left and the right integrals, defined on the subalgebra $A_{d0}$. They also find that the modular element, relating the left with the right integral, is equal to (a scalar multiple of) $F^2$ where the element $F$ in the multiplier algebra of the subalgebra that was already constructed within the theory of compact matrix quantum groups, using the orthogonality relations of the matrix elements of the irreducible corepresentations. 
\ssnl
Finally in Theorem 3.4, they essentially show that the dual of $A_d$ gives again $A_c$.
\ssnl
So what is done by Podle\'s and Woronowicz in these two sections of \cite{Po-Wo} can be summarized as follows. They define a discrete quantum group as a pair of a C$^*$-algebra $A$ with a coproduct $\Delta$ on $A$, where $A$ is a direct sum of matrix algebras and so that furthermore there is a counit and an antipode. They associate to each compact quantum matrix group a dual discrete quantum group and they obtain further results about this discrete quantum group (such as the existence of the integrals) using this fact and various results about compact quantum groups and their corepresentation theory.
\nl
\bf Comparison with the multiplier Hopf algebra approach  \rm 
\nl
In \cite{VD-discrete}, a discrete quantum group is defined as a multiplier Hopf $^*$-algebra $(A,\Delta)$ where $A$ is a (algebraic) direct sum of matrix algebras (see Definition 2.3 in \cite{VD-discrete}). A multiplier Hopf $^*$-algebra is defined in \cite{VD-mha}, Definitions 2.3 and 2.4.
\ssnl
This definition is not very different from the one given by Podle\'s and Woronowiz in the sense that it takes little effort to pass from one to the other. Indeed, if we start e.g.\ with a multiplier Hopf $^*$-algebra, the general theory allows to construct a counit and an antipode (see Section 3 and 4 in \cite{VD-mha}). Then, as the algebra is a direct sum of matrix algebras, it is not hard to extend the coproduct and the counit to the C$^*$-direct sum. This takes us to a discrete quantum group as introduced (implicitly) in Theorem 3.1 of \cite {Po-Wo}.
\ssnl
Conversely, if we start with a discrete quantum group as in \cite{Po-Wo},  we simply restrict the coproduct to the algebraic direct sum of the components and we obtain a discrete quantum group in the sense of \cite{VD-discrete}. The counit and the antipode are used to show that the basic axioms for a multiplier Hopf $^*$-algebra (see \cite{VD-mha} and/or \cite{VD-discrete}) are fulfilled. It is worthwhile noticing that the axioms for multiplier Hopf algebras can also be formulated in terms of the existence of a counit and an antipode (as for Hopf algebras) and then the two definitions come even closer together.
\ssnl
So the difference between the two approaches does not really lie in the definition of what is a discrete quantum group, but rather in other aspects. We will discuss these now.
\ssnl
The first, probably the main difference is that in \cite{VD-discrete} the existence of the invariant integrals is proven from the axioms. In \cite{Po-Wo}, the existence of the integrals on the discrete quantum group is obtained given the fact that it is the dual of a compact quantum group. However, if you want to show that every discrete quantum group is the dual of a compact one, you need the integrals to do that.
\ssnl
Along the same lines, by treating the discrete quantum groups ad duals of compact quantum groups, also properties of discrete quantum groups are proven using results about compact quantum groups. In the other, independent approach, all results about discrete quantum groups are obtained from the axioms.
\ssnl
The second difference, of another nature but maybe equally important, is the fact that we consider the duality between compact and discrete quantum groups as a special case of the duality for algebraic quantum groups. In this theory, the general case turns out to be simpler than the special case. If e.g.\ the general construction of the dual (as for algebraic quantum groups) is applied to the special cases of a compact or a discrete quantum group, it is easy to see that the result is a discrete, respectively a compact quantum group. The arguments are simpler than the ones given in \cite{Po-Wo}. 
\ssnl
This leads in general, also at other places where recently other aspects of compact and discrete quantum groups are studied, to simplifications and a better understanding.
This should not come as a surprise because it is already the case in the approach to discrete quantum groups itself as treated in  \cite{Po-Wo}.
\nl
\bf Compact quantum groups as duals of discrete quantum groups \rm 
\nl
As one can study discrete quantum groups as duals of compact quantum groups, similarly it is possible to study compact quantum groups as duals of discrete quantum groups. This is possible because on the one hand, the discrete quantum groups, as well as the compact quantum groups can be studied on their own, while on the other hand, it is shown in each of the two cases that the dual of one type is an object of the other type. However, it must be said in order to get the full picture, one needs to develop each of the theories, up to a certain extend, separately. One has e.g.\ to prove the existence of the integrals in each of the two cases.
\ssnl
Starting with a discrete quantum group and from there, studying compact quantum groups is indeed possible (and is somewhat easier than the other, more common way). We show this in these notes. But this only makes sense if you know already that any compact quantum group is the dual of a discrete quantum group. And as mentioned already, this result requires to develop the compact quantum groups independently. But the same is true for the other direction.
\snl
Taken abstraction of this (and of the historical context), people who want to work with compact and discrete quantum groups can consider any of the two approaches. To a great extend, it is a matter of taste which you take as a starting point (although we are convinced that starting with discrete quantum groups is a bit easier). However, for a full understanding, we feel it is worthwhile to  study the two theories also as independent entities.

\nl\nl





\bibliographystyle{amsplain}

\bibliography{references.bib}

\end{document}